\documentclass[11pt]{article}
\usepackage{amssymb}
\usepackage[top=23mm,bottom=23mm,left=23mm,right=23mm]{geometry}
\title{Lifschitz singularity for subordinate Brownian motions in presence of the Poissonian potential on the Sierpi\'nski gasket}
\author{ Kamil Kaleta
 and  Katarzyna Pietruska-Pa\l{}uba}

\usepackage{color}

\begin{document}
\maketitle

\begin{abstract}
We establish the Lifschitz-type singularity around the bottom of the spectrum for the integrated density of states for a class of subordinate
Brownian motions in presence of the nonnegative Poissonian random potentials, possibly of infinite range, on the Sierpi\'nski gasket. We also study the long-time behaviour for the corresponding averaged Feynman-Kac functionals.

\medskip
\noindent
2010 {\it MS Classification}: Primary 60J75, 60H25, 60J35; Secondary 47D08, 28A80.

\end{abstract}

\footnotetext{K. Kaleta \\ Institute of Mathematics, University of Warsaw, ul. Banacha 2, 02-097 Warszawa and Institute of Mathematics and Computer Sciences, Wroc{\l}aw University of Technology, Wyb. Wyspia\'nskiego 27, 50-370 Wroc{\l}aw, Poland \\ \emph{ e-mail: kkaleta@mimuw.edu.pl, kamil.kaleta@pwr.wroc.pl}}

\footnotetext{K. Pietruska-Pa{\l}uba \\ Institute of Mathematics, University of Warsaw, ul. Banacha 2, 02-097 Warszawa, Poland \\
\emph{ e-mail: kpp@mimuw.edu.pl}}

\footnotetext{
K. Kaleta was supported by the National Science Center (Poland) internship
grant on the basis of the decision No. DEC-2012/04/S/ST1/00093 and by the Foundation for Polish Science.
}

\newcommand{\rd}{\mathbb{R}^d}
\newcommand{\dist}{\emph{dist}}
\newcommand{\tr}{\mbox{tr}}
\newtheorem{theo}{\bf Theorem}[section]
\newtheorem{coro}{\bf Corollary}[section]
\newtheorem{lem}{\bf Lemma}[section]
\newtheorem{rem}{\bf Remark}[section]
\newtheorem{defi}{\bf Definition}[section]
\newtheorem{exam}{\bf Example}[section]
\newtheorem{fact}{\bf Fact}[section]
\newtheorem{prop}{\bf Proposition}[section]
\newtheorem{oq}{\bf Open question}

\providecommand{\pro}[1]{(#1_t)_{t \geq 0}}
\providecommand{\proo}[1]{(#1_t)_{t \in \R}}
\providecommand{\semi}[1]{\{#1_t: t \geq 0\}}

\providecommand{\bri}[1]{(#1_r)_{0 \leq r < t}}
\providecommand{\seq}[1]{(#1_n)_{n\in \mathbb{N}}}

\newcommand{\cS}{\mathcal{S}}
\newcommand{\cmw}{\mathcal{w}}
\newcommand{\cms}{\mathcal{s}}
\newcommand{\cQ}{\mathcal{Q}}
\newcommand{\cG}{\mathcal{G}}
\newcommand{\cD}{\mathcal{D}}
\newcommand{\cA}{\mathcal{A}}
\newcommand{\cL}{\mathcal{L}}
\newcommand{\cK}{\mathcal{K}}
\newcommand{\cH}{\mathcal{H}}
\newcommand{\cF}{\mathcal{F}}
\newcommand{\cE}{\mathcal{E}}
\newcommand{\cB}{\mathcal{B}}
\newcommand{\cX}{\mathcal{X}}
\newcommand{\cY}{\mathcal{Y}}
\newcommand{\cT}{\mathcal{T}}
\newcommand{\cW}{\mathcal{W}}
\newcommand{\cM}{\mathcal{M}}
\newcommand{\cN}{\mathcal{N}}
\newcommand{\pr}{\mathbf{P}}
\newcommand{\ex}{\mathbf{E}}
\newcommand{\qpr}{\mathbb{Q}}
\newcommand{\qex}{\mathbb{E}_\mathbb{Q}}
\newcommand{\R}{\mathbf{R}}
\newcommand{\1}{\mathbf{1}}
\newcommand{\Z}{\mathbf{Z}}
\newcommand{\trmu}{\mathbf{\mu}}
\newcommand{\loc}{\footnotesize{loc}}
\newcommand{\Rd}{\mathbf{R}^d}
\newcommand{\N}{\mathbf{N}}

\newcommand{\p}{\phantom{aa}}

\makeatletter \@addtoreset{equation}{section}
\renewcommand{\theequation}{\thesection.\arabic{equation}}
\makeatother
 \maketitle

\section{Introduction}
The integrated density of states is one of central objects in
the physics of large-volume systems, especially systems with
in-built randomness. The randomness can come from the interaction
 with external force field, described by its potential $V.$ This leads us to the study of random Hamiltonians, in particular those of Schr\"{o}dinger type: given a sufficiently regular, possibly random, potential $V$ one considers the operator
\[H:= H_0 + V,\]
where $H_0$ is the Hamiltonian of the system with no potential interaction. The best analyzed situation is that of $H_0=-\Delta$
(in various state-spaces $X$). The spectrum of $H$ is typically not discrete.
Moreover, spectral properties of such infinite-volume (i.e. defined with the whole space $X$) Schr\"odinger operators are
usually difficult to handle. The notion of the integrated density of states can come to the rescue: it captures some of the properties
of the spectral distribution,
while being easier to calculate and
easier to work with \cite[Chapter VI]{bib:Car-Lac}.

Informally speaking, one considers operators $H$ restricted to
a finite volume $\Omega\subset X,$ build empirical measures $l_\Omega$ based on the spectra of these operators normalized by the volume of $\Omega$, and then one takes the limit of $l_\Omega,$ in appropriate sense, when $\Omega\nearrow X.$ The resulting limit (if it exists) is called the integrated density of states (IDS, for short). Same procedure can be performed for random potentials $V^\omega$ -- in this case one is interested in the almost-sure
limiting behavior of measures $l_\Omega.$ When the potential $V^\omega$ exhibits some ergodicity properties then the limit can be nonrandom.

This paper is concerned with random Schr\"odinger operators with nonnegative Poissonian potentials. In this case, the existence of the nonrandom IDS is a common feature and for $H_0=-\Delta$ has been proven e.g. in
the Euclidean space \cite{bib:Nakao}, hyperbolic space \cite{bib:Szn-hyp1}, the Sierpi\'{n}ski gasket \cite{bib:KPP-PTRF},
other nested fractals \cite{bib:Sh}. In all these situations one has
the so-called Lifschitz singularity: the rate of decay of the IDS at the bottom of the spectrum is faster than this of the IDS for the system without external interaction. Note also that the Lifschitz singularity is closely related to the behaviour of the so-called Wiener sausage when $t\to\infty$ (for the sausage asymptotics in the classical case see \cite{bib:DV}).

While the IDS based on the Laplacian is fairly well understood
(see e.g. \cite{bib:Car-Lac}, \cite{bib:Stol}), it is not so
for the IDS based on nonlocal operators. In the case of L\'{e}vy processes on $\mathbb R^d,$ the existence and asymptotical properties of IDS with Poissonian potentials have been established in \cite{bib:Okura,bib:Okura81}. Up to date, there were no results concerning the `nonlocal IDS' on irregular sets, such as fractals. Recently, we have proven the existence of the IDS for subordinate Brownian motions on the Sierpi\'nski gasket perturbed by Poissonian potentials with two-argument profiles $W$ that may have infinite range and local singularities \cite{bib:KaPP}. The Lifschitz tail for stable processes on the Sierpi\'nski gasket evolving among killing Poissonian obstacles was derived in \cite{bib:Kow-KPP}. The present paper is meant as the continuation
of \cite{bib:KaPP} in the potential case: under appropriate assumptions on the potential $V$
(expressed in terms of its profile function $W$) and the Laplace exponent of the subordinator $S$ (assumed to be a complete Bernstein function), we analyse the asymptotical behaviour of the
IDS based on the generator of the resulting subordinate Brownian motion evolving in presence of the potential $V$.
We first estimate the Laplace transform of the IDS (Theorems \ref{th:lower} and \ref{th:upper-gen}), and then
use exponential Tauberian theorems from \cite{bib:F} to trasform
them into estimates on the IDS itself (Theorems \ref{th:lif-upper} and \ref{th:lif-lower}). The proof of lower bound seems to be easier,
while for the proof of upper bound we need to reduce the problem to study the subordinate Brownian motion reflected in a gasket triangle of size $2^M$ perturbed by some special periodization of the initial potential (Lemma \ref{lem:monotone}), and compare the principal eigenvalue of its generator with the principal eigenvalue related to the stable process reflected in the unit triangle, with  rescaled potential (Lemmas \ref{lem:reflected_dir} and \ref{lem:upper_basic}). After this simplification, we can just proceed with stable process and employ the coarse-graining technique (`the enlargement of obstacles method') of Sznitman \cite{bib:Szn1}, adapted recently to the case of non-diffusive processes in \cite{bib:Kow}. In present work, Sznitman's theorem is needed for the potential case (the original work \cite{bib:Szn1} was concerned with killing obstacles).
To make the paper self-contained, we give a proof of the desired theorem (Theorem \ref{compare}) in the Appendix. We also derive the estimates for the corresponding averaged Feynman-Kac functional, which can be interpreted as the survival probability  of the process killed by the potential $V$ up to time $t$ (Theorem \ref{th:lower-funct-1} and Remark \ref{re:func_upper}).
Our proof of the upper bound hinges on properties of the reflected subordinate Brownian motions on the gasket. Construction of such processes relies on the specific geometry of the gasket and does not seem to have an obvious extension to more general fractals.

A remarkable feature of our results is that they take into account also the long range interaction which comes from the decay rate of the profile $W$ at infinity and often has a decisive impact on the properties of the IDS. This seems to be new even in the case of Brownian motion on the Sierpi\'nski gasket. In the jump case, our bounds reflect very well the competition between intensity of large jumps of the process and the rate of killing of the Poissonian potential given by the tail properties of the potential profile. Our approach covers a wide range of jump subordinators with drift (the resulting subordinate process is then often called \emph{jump diffusion}), as well as purely jump subordinators including stable, mixture of stable, some logarithmic-stable and others (see the examples in Section \ref{sec:examples}). Unfortunately, our Theorem \ref{th:upper-gen} cannot be applied to the relativistic stable subordinators. We believe that this result holds true in this case as well, but the proof would require tools specialized to those specific processes, not available yet.

\section{Preliminaries}

\subsection{The Sierpi\'nski gasket} \label{subsec:gasket}
The infinite Sierpi\'nski Gasket we will be working with is
defined as a blowup of the unit gasket, which in turn is defined
as the fixed point of the hyperbolic iterated function system
in $\mathbb R^2,$ consisting of three maps:
\[
\phantom{}\ \ \ \ \ \ \phi_1(x)=\frac{x}{2},\ \ \ \ \   \phi_2(x)=
\frac{x}{2}+\left(\frac{1}{2},0\right) \ \ \ \ \ \
\phi_3(x)=\frac{x}{2}+\left(\frac{1}{2}, \frac{\sqrt 3}{2}\right). \hfill \]
The unit gasket, $\mathcal G_0,$ is the unique compact subset of
$\mathbb{R}^2$ such that
\(\mathcal G_0= \phi_1(\mathcal G_0)\cup \phi_2(\mathcal G_0)\cup\phi_3(\mathcal G_0).\) Then
we set:
\[\mathcal G_n=2^n \mathcal G_0 =((\phi_1^{-1}))^n(\mathcal G_0),\;\;\;\;\;\mathcal G=\bigcup_{n=1}^\infty \mathcal G_n.\]
All the  triangles of
size $2^M,$ $ M\in\mathbb Z, $ that build up the infinite gasket will be denoted by $\mathcal T_M,$ and the collection of their vertices -- by
 $\mathcal V_M.$

We equip the gasket with the shortest path distance
$d(\cdot,\cdot)$: for $x,y\in \bigcup_M\mathcal V_M,$ $d(x,y)$ is the infimum of
Euclidean lengths of all paths, joining $x$ and $y$ on the gasket.
For general $x,y\in \cG,$  $d(x,y)$ is obtained by a limit
procedure. This metric is equivalent to the usual Euclidean metric
inherited from the plane.
    Observe that
$\mathcal G_M= B(0, 2^M),$ where the ball is taken in either the
Euclidean or the shortest path metric.

By $m$ we denote the Hausdorff measure on $\mathcal G$ in
dimension $d=\frac{\log 3}{\log 2},$   normalized to have
$m(\mathcal G_0)=1.$ The number $d$ is called the fractal
dimension of $\mathcal G.$ Another characteristic number of
$\mathcal G,$ namely $d_w=\frac{\log 5}{\log 2}$ is called the
walk dimension of $\mathcal G.$ The spectral dimension of
$\mathcal G$ is $d_s=\frac{2d}{d_w}.$
The measure $m$ is a $d-$measure, i.e. there exists  a constant $c_{2.1}>0$ such that for all $x\in\mathcal G,$ $r>0$
one has
\begin{equation}\label{eq:d-meas}
c_{2.1}^{-1} r^d\leq m(B(x,r)) \leq c_{2.1} r^d
\end{equation}

By an elementary slicing argument we obtain the following bound, valid for $x\in\mathcal G, $ $ \lambda>0,$ $a>0:$
\begin{equation}\label{eq:meas-int}
c_{2.2}^{-1} \frac{1}{a^\lambda}\leq \int_{d(x,y)>a} \frac{{\rm d}m(y)}{d(x,y)^{d+\lambda}} \leq c_{2.2} \frac{1}{a^\lambda},
\end{equation}
where $ c_{2.2}=c_{2.2}(d,\lambda)>0$ is a numerical constant.

In the sequel, we will need a projection from $\mathcal G$ onto
$\mathcal G_M,$ $M=0,1,2,...$ To define it properly, we first put
labels on the set $\mathcal V_0$ (see \cite{bib:KPP-PTRF}).

Observe that $\mathcal
G_0\subset(\mathbb{Z}_+)e_1+(\mathbb{Z}_+)e_2,$ where $e_1=(1,0)$
and $e_2=\left(\frac{1}{2},\frac{\sqrt 3}{2}\right).$ Next,
consider the commutative $3-$group $\mathbb{A}_3$ of even
permutations of 3 elements,  $\{a,b,c\},$ i.e. $\mathbb{A}_3=\{id,
(a,b,c), (a,c,b) \},$ and we denote $p_1=(a,b,c),$ $p_2=(a,c,b).$
The mapping
\[\mathcal V_0\ni x=ne_1+me_2\mapsto p_1^n\circ p_2^m\in\mathbb{A}_3\]
is well defined, and for $x\in \mathcal V_0$ we put $l(x)= (p_1^n\circ
p_2^m)(a).$ This way, every triangle of size 1 has its vertices
labeled `$a,b,c$', see   \cite[Fig. 4 and 5]{bib:KPP-PTRF}. Note that this property extends to every
triangle of size $2^M,$ which corresponds to putting labels on the
elements of $\mathcal V_M:$ every triangle of size $2^M$ has three distinct labels on its vertices.

Let $M\geq 0 $ be fixed. For $x\in \mathcal G\setminus \mathcal V_M,$ there
is a unique triangle of size $2^M$ that contains $x,$
$\Delta_M(x),$ and so $x$ can be written as $x=x_a a(x)+x_b
b(x)+x_c c(x),$ where $a(x),$ $b(x),$ $c(x)$ are the vertices of
$\Delta_M(x)$ with labels $a,b,c$  and $x_a,x_b,x_c\in (0,1),$
$x_a+x_b+x_c=1.$ Then we define the projection:
\[\mathcal G\setminus \mathcal V_M\ni x\mapsto \pi_M(x)= x_a\cdot a(M)+x_b\cdot b(M) +x_c\cdot c(M)\in \mathcal G_M, \]
where $a(M),$ $b(M),$
 $c(M)$ are the vertices of the triangle $\mathcal G_M$ with
 corresponding labels $a,b,c.$
 When $x\in \mathcal V_M,$ then $x$ itself has a label assigned and it
 then mapped onto the corresponding vertex of $\mathcal G_M.$

\subsection{Subordinate Brownian motions and their Schr\"odinger perturbations}

Let $\mathcal G^*$ be the two-sided infinite gasket, i.e. the set $\mathcal G^*=\mathcal G\cup i(\mathcal G),$ where $i$ is the reflection
of $\mathbb R^2$  with respect to the $y-$axis. Denote by $\tilde Z=(\tilde Z_t, \pr_x)_{t\geq 0,x\in\mathcal G^*}$ the Brownian motion on
$\mathcal G^*$, as defined in \cite{bib:BP}. It is a strong Markov and Feller process, whose transition density with respect to the Hausdorff measure is symmetric in its space variables, continuous, and fulfils
 the following subgaussian estimates:
\begin{eqnarray}
\label{eq:gaussian}
c_{2.3} t^{-d_s/2} {\rm e}^{-c_{2.4} \left(\frac{d(x,y)}{t^{1/d_w}}\right)^{d_w/(d_w-1)}} \leq \tilde g(t,x,y) \leq c_{2.5} t^{-d_s/2} {\rm e}^{-c_{2.6} \left(\frac{d(x,y)}{t^{1/d_w}}\right)^{d_w/(d_w-1)}}\hskip -1cm,\\
 \nonumber \phantom{} \hskip 7cm \quad x,y \in \mathcal G^*, \quad t>0,
\end{eqnarray}
with positive constants $c_{2.3},...,c_{2.6}.$
By $Z$ we denote the Brownian motion on the one-sided Sierpi\'{n}ski gasket $\mathcal G$, which is obtained from $\tilde Z$ by the projection $\mathcal G^*\to \mathcal G$. One can directly check that its transition densities are given by
\(g(t,x,y)= \tilde g(t,x,y)+\tilde g(t,x,i(y))\) for $x\neq 0$ and twice this quantity when $x=0.$ The functions $g$ share most of the properties of $\tilde g$, including continuity and the subgaussian estimates, with possibly worse constants $c_{2.3}-c_{2.6}$. We stick to the estimate (\ref{eq:gaussian}) for $g(t,x,y)$
as well.

Let $S=(S_t, \pr)_{t \geq 0}$ be a subordinator, i.e. an increasing L\'{e}vy process taking values in $[0,\infty]$ with $S_0=0$.
Denote $\eta_t({\rm d}u) = \pr(S_t \in du)$, $t \geq 0$. As usual, if the measures $\eta_t({\rm d}u)$ are absolutely continuous with respect to the Lebesgue measure, then the corresponding densities are also denoted by $\eta_t(u)$. The law of the subordinator $S$ is uniquely determined by its Laplace transform $\int_0^{\infty} {\rm e}^{-\lambda s} \eta_t({\rm d}s) ={\rm e}^{-t\phi(\lambda)},$ $\lambda>0.$
The function $\phi:(0,\infty) \to [0,\infty)$ is called the Laplace exponent of $S$ and can be represented as
\begin{equation}
\label{def:phi}
\phi(\lambda) = b \lambda + \psi(\lambda) \quad \mbox{with} \quad \psi(\lambda)=\int_0^{\infty} (1-{\rm e}^{-\lambda s})\rho({\rm d} s),
\end{equation}
where $b\geq 0$ is called the drift term and $\rho$, called the L\'evy measure of $S$, is a $\sigma$-finite measure on $(0,\infty)$ satisfying $\int_0^{\infty} (s \wedge 1) \rho({\rm d}s) < \infty$. It is well known that $\phi$ is  a Bernstein function. For more details of subordinators and Bernstein functions we refer the reader to \cite{bib:Ber, bib:SSV}.

We always assume that $Z$ and $S$ are independent. The process $X=(X_t, \pr_x)_{t\geq 0,x\in\mathcal G}$ given by
$$
X_t := Z_{S_t}, \quad t \geq 0,
$$
is called the subordinate Brownian motion on $\cG$ (via subordinator $S$). It is a symmetric Markov process having c\`adl\`ag paths. Its natural filtration is always assumed to fulfil the usual conditions.

For the entire paper, we make the following  assumptions \textbf{(S1)} and \textbf{(S2)} on the subordinator.
\begin{itemize}
\item[\textbf{(S1)}] $\forall\;t<0$ one has
$
\eta_t(\left\{0\right\})=0 $ and $ \int_{0^{+}}^\infty \frac{1}{u^{d_s/2}}\,\mathbb \eta_t({\rm d}u)=:c_{2.7}(t) <\infty, \nonumber
$
\item[\textbf{(S2)}] $\forall\;t<0$ one has
$
\int_1^{\infty} \eta_t(u,\infty) \frac{{\rm d} u}{u}<\infty.
$
\end{itemize}

Assumptions \textbf{(S1)}--\textbf{(S2)} are satisfied by a wide class of subordinators. Observe that whenever $S_t=ct$ and $X_t=Z_{ct}$ for some $c>0$, i.e., $X$ is the time-rescaled Brownian motion with $\eta_t({\rm d}u)=\delta_{ct}({\rm d}u)$ and $\phi(\lambda)=c\lambda$, they immediately hold. Some other specific examples including jump processes and sufficient conditions for them are discussed in \cite[Remark 2.1, Lemma 2.2 and Example 2.1]{bib:KaPP}. See also Section \ref{sec:examples} at the end of this paper.  For further examples we refer the reader to \cite{bib:SSV}.

In the sequel we will need the following estimate of the tail of the subordinator $S$ at infinity. It can also be used to verify \textbf{(S2)} for given  $\phi.$

\begin{lem} \label{lem:tails}
Let $S$ be a subordinator with Laplace exponent $\phi$  such that $\mathbf P[S_t=0]=0,$ for every $t>0$. Then we have
$$
\eta_t(A, \infty) \leq \frac{t}{1-e^{-1}} \, \int_0^{\frac{1}{A}} \frac{\phi(\lambda)}{\lambda} {\rm d} \lambda, \quad t>0, \ A >0.
$$
\end{lem}
\noindent{\bf Proof.}
Standard arguments as in \cite[Lemma 2.2]{bib:KaPP}, yield
$$
\int_0^{\infty} e^{-\lambda u} \eta_t(u, \infty) {\rm d}u \leq \frac{t \phi(\lambda)}{\lambda}, \quad \lambda >0, \ t>0,
$$
which, by monotonicity, leads to
$$
e^{-\lambda A} \, A \, \eta_t(A, \infty) {\rm d}u \leq \int_0^{A} e^{-\lambda u} \eta_t(u, \infty) {\rm d}u \leq \frac{t \phi(\lambda)}{\lambda}, \quad \lambda >0, \ t>0, \, A>0.
$$
By integrating  both sides of the above inequality with respect to $\lambda$ over $(0,1/A)$, we thus get
$$
(1-e^{-1}) \eta_t(A, \infty) \leq t \int_0^{\frac{1}{A}} \frac{\phi(\lambda)}{\lambda} {\rm d} \lambda, \quad \lambda >0, \ t>0, \, A>0,
$$
which is the claimed inequality. \hfill$\Box$

\smallskip

By the first part of \textbf{(S1)}, the process $X$ has symmetric and strictly positive transition densities given by
\begin{equation}\label{eq:subord}
p(t,x,y) = \int_0^{\infty} g(u,x,y) \eta_t({\rm d}u), \quad t>0, \quad x, y \in \cG,
\end{equation}
while the second part guarantees that
\begin{equation}\label{eq:sup-of-g}
\sup_{x,y \in \cG} p(t,x,y) \leq c_{2.8} (t)< \infty \quad \mbox{for every} \quad t >0 \quad \mbox{and} \quad c_{2.9}:=\sup_{t \geq 1} c_{2.8}(t)< \infty.
\end{equation}
Moreover, for each fixed $t >0$, $p(t,\cdot,\cdot)$ is a continuous function on $\cG \times \cG$, and for each fixed $x, y \in \cG$, $p(\cdot,x,y)$ is a continuous function on $(0,\infty)$. The general theory of subordination (see, e.g., \cite[Chapter 12]{bib:SSV}) yields that the process $X$ is a Feller process and, in consequence, a strong Markov process. Also, by \textbf{(S1)}, it has the strong Feller property.

For an open set $U \subset \cG$ by $\tau_U:=\inf\left\{t \geq 0: X_t \notin U \right\}$ we denote the first exit time of the process $X$ from $U$, while $T_{D}=\tau_{D^c}$ denotes the first entrance time into the closed set $D$. By $(P_t)_{t\geq 0}$ we denote the semigroup with kernel $p(t,x,y),$ and by $(P_t^U)_{t\geq 0}$ -- the semigroup related to the process killed \linebreak on exiting an open set $U$; $\lambda_1(U)$ is the principal eigenvalue of its generator.
Finally, by $\pr^t_{x,y}$ we \linebreak denote the bridge measures corresponding to process $X$ on $D([0,t],\cG)$ (for more details we refer to \cite[p. 9]{bib:KaPP}).

We say that a Borel function $V:\mathcal G\to\mathbb R$ is in Kato class $\cK^X$ related to the process $X$ if
\begin{equation}
\label{eq:Kato}
\lim_{t \searrow 0} \sup_{x \in \cG} \int_0^t \ex_x |V(X_s)|\,{\rm d}s = 0.
\end{equation}
Also, $V \in \cK^X_{\loc}$ (local Kato class), when $\1_B V \in \cK^X$ for every ball $B \subset \cG$. One can show that $L^{\infty}_{\loc}(\cG) \subset \cK^X_{\loc} \subset L^1_{\loc}(\cG,m)$.

In this paper, we study the subordinate Brownian motions on $\cG$ perturbed by random Schr\"odinger potentials of  the Poissonian type, which are defined by
\begin{equation}
\label{eq:poiss0} V(x,\omega):= \int_{\cG} W(x,y)
\mu^{\omega}({\rm d}y), \quad x \in \cG, \ \ \omega \in \Omega,
\end{equation}
where $\mu^\omega$ is the random counting measure corresponding to
the Poisson point process on $\cG$,  with intensity $\nu {\rm d}m$,  $\nu>0,$ defined on some
probability space $(\Omega,\mathcal M, \mathbb{Q})$,  and $W:\cG \times \cG \to \mathbb R_+$ is a
measurable, nonnegative profile function.
We note for later use that for any measurable function $f:\cG\to \mathbb R_+$ one has
\begin{equation}\label{eq:exp-formula}
\mathbb{E_Q}\left[{\rm e}^{-\int_\cG f(y)\,\mu^{\omega}({\rm d}y)}\right] = {\rm e}^{-\int_{\cG} \nu (1-{\rm e}^{-f(y)})m({\rm d} y)}.
\end{equation}Throughout the paper we assume that the Poisson process and the Markov process $X$ are independent and impose the following regularity assumptions on the profile function.

\smallskip

\begin{itemize}
\item[\textbf{(W1)}]
$W\geq 0$, $W(\cdot,y) \in \cK_{\loc}^X$ for every $y \in \cG$ and there exists a function $h \in L^1(\cG,m)$ such that $W(x,y) \leq h(y)$, whenever $d(y,0)\geq 2d(x,0)$.
\item[\textbf{(W2)}] $\sum_{M=1}^{\infty} \sup_{x \in \cG} \int_{B(x,2^{M/4})^c} W(x,y) dm(y) < \infty$
\item[\textbf{(W3)}] there is $M_0 \in \mathbb Z$ such that
\begin{equation}\label{eq:wu}
\sum_{y{'} \in \pi_M^{-1}(\pi_M(y))} W(\pi_M(x),y{'}) \leq \sum_{y{'}
\in \pi_M^{-1}(\pi_M(y))} W(\pi_{M+1}(x),y{'}), \quad x, y \in \cG,
\end{equation}
for every $M \in \mathbb Z$, $M \geq M_0$.
\end{itemize}
The conditions \textbf{(W1)}--\textbf{(W3)} have been recently introduced in \cite{bib:KaPP}. Under \textbf{(W1)} we have that $V(\cdot,\omega) \in \cK^X_{\loc}$ for $\qpr$-almost all $\omega \in \Omega$, while the remaining conditions \textbf{(W2)}--\textbf{(W3)} guarantee  sufficient regularity, needed to study the spectral problem that we address in the present paper. For discussion of the above assumptions we refer the reader to \cite[Subsection 2.3.2]{bib:KaPP}.

For the Poissonian random potential $V$, we consider the Feynman-Kac semigroups $(T^{D,M,\omega}_t)_{t\geq 0}$ related to the process killed outside $\cG_M$, $M \in \mathbb Z$, consisting of the operators:
\begin{equation}
\label{def:sem-dir}
T_t^{D,M,\omega} f(x) = \ex_x \left[{\rm e}^{-\int_0^t V(X_s, \omega) {\rm d}s} f(X_t); t<\tau_{\cG_M}\right], \quad f \in L^2(\cG_M,m), \quad M \in \mathbb Z, \quad t>0.
\end{equation}
For $\qpr$-almost all $\omega \in \Omega$ and every $t>0$, $T_t^{D,M,\omega}$ are symmetric, ultracontractive and Hilbert-Schmidt operators admitting  measurable, symmetric and bounded kernels $u^M_D(t,x,y)$ which are known to have the following very useful bridge representation \cite[(2.27)]{bib:KaPP}
\begin{equation}
\label{def:sem-dir-kernel-bridge}
u^M_D(t,x,y) = p(t,x,y) \ \ex^t_{x,y}\left[{\rm e}^{-\int_0^t V(X_s,\omega) {\rm d}s} ; t<\tau_{\cG_M}\right],  \quad M \in \mathbb Z, \quad x,y \in \cG_M, \quad t>0.
\end{equation}
Denote by $A^{D,M,\omega}$ the $L^2(\cG,m)$--generator of the semigroup $(T^{D,M,\omega}_t)_{t\geq 0}$.
By analogy to the Euclidean case, the operator $-A^{D,M,\omega}$ is called the \emph{generalized Schr\"odinger operator} corresponding to the generator of the process $X$ with Dirichlet (outer) conditions. The spectrum of $-A^{D,M,\omega}$ is purely discrete. The corresponding eigenvalues can be ordered as $0\leq \lambda_1^{D,M}(\omega)\leq \lambda_2^{D,M}(\omega)\leq ... \to \infty$. For discussion and verification of the properties and facts listed above we refer the reader to \cite[Subsection 2.3.1 and 2.3.2]{bib:KaPP}.

The basic objects we consider are the random empirical measures on $\mathbb R_+:=[0,\infty)$ based on the spectra of $-A^{D,M,\omega}$, normalized by the volume of $\cG_M$:
\begin{equation}\label{eq:el-d}
l_M^D(\omega):= \frac{1}{m(\cG_M)} \sum_{n=1}^{\infty}
\delta_{\lambda_n^{D,M}(\omega)}, \quad M \in \mathbb Z_+,
\end{equation}
and their Laplace transforms $L_M^D(t,\omega) := \int_0^{\infty} {\rm e}^{-\lambda t} {\rm d}l_M^D(\omega)(t) = \frac{1}{m(\cG_M)} \sum_{n=1}^{\infty} {\rm e}^{-\lambda_n^{D,M}(\omega) t}$ which have the following representation
\begin{eqnarray*}
L_M^D(t,\omega)=
\frac{1}{m(\cG_M)} \tr T_t^{D,M,\omega}
= \frac{1}{m(\cG_M)} \int_{\cG_M} p(t,x,x) \ex^t_{x,x}\left[e^{-\int_0^t V(X_s,\omega) ds}; t < \tau_{\cG_M} \right] m({\rm d} x).
\end{eqnarray*}
We have recently proven the following result on the convergence of $L_M^D$ and $l_M^D$ as $M \to \infty$.

\begin{theo}\label{th:meanlimit} {\rm \cite[Theorems 3.1 and 3.2]{bib:KaPP}}
Let $S$ be a subordinator satisfying the assumptions \textbf{(S1)}-\textbf{(S2)} and let $V$ be a Poissonian random field with the profile $W$ satisfying the conditions \textbf{(W1)}-\textbf{(W3)}. Then
for every $t>0$, $\mathbb{E}_\qpr[L^D_M(t,\omega)]$ converges to a finite limit $L(t)$ as $M\to\infty.$ Moreover, $\qpr-$almost surely, the random measures $l_M^D(\omega)$
 vaguely converge to a common nonrandom limit measure $l$ on $\mathbb R_+,$ with Laplace transform $L(t)$.
\end{theo}
The deterministic measure $l$ given by the above theorem is called the \emph{integrated density of states} (IDS) for the process $X$ perturbed by the Poissonian potential $V$ on $\cG$. The present paper is devoted to study of the limiting behaviour of $l[0,\lambda)$ when $\lambda \to 0^{+}$.

As stated above, the assumptions \textbf{(S1)}--\textbf{(S2)} and \textbf{(W1)}--\textbf{(W3)} guarantee the existence of IDS in our settings. In  present paper, the conditions \textbf{(S1)}--\textbf{(S2)} only give the general framework for our study. From now on, we will restrict our attention to the class of the so-called \emph{complete} subordinators and impose some additional regularity on the corresponding Laplace exponent $\phi$. Recall that the subordinator $S$ is called \emph{complete} if its Laplace exponents $\phi$ is \emph{complete Bernstein function}, i.e., the corresponding L\'evy measure $\nu$ is absolutely continuous with respect to the Lebesgue measure with completely monotone density (see e.g. \cite[Chapter 6]{bib:SSV}).

\section{Lower bounds}\label{seclower}

\subsection{Lower bounds for the integrated density of states}

As explained above, the integrated density of states $l$ is the vague limit of the empirical measures based on the Laplacians
on $\mathcal G_M$ with Dirichlet boundary conditions, and its Laplace transform $L(t)$  for any given $t>0$ can be expressed as the limit:
\begin{equation} \label{eq:L-lim}
L(t)=\lim_{M\to\infty}\mathbb{E_Q}[L^D_M(t,\omega)],
\end{equation}
where
\[L^D_M(t,\omega)=\frac{1}{m(\mathcal G_M)}\int_{\mathcal G_M} p(t,x,x) \mathbf E_{x,x}^t\left[ {\rm e}^{-\int_0^t V(X_s,\omega)\,{\rm d} s} \mathbf{1}_{\{\tau_{\mathcal G_M}>t\}}\right]\,{\rm d} m(x).\]

Before we proceed, we introduce some notation. For a gasket triangle $\Delta,$ let us define
\begin{equation}\label{eq:L-delta}
L^\Delta(t,\omega)=\frac{1}{m(\Delta)}\int_{\Delta} p(t,x,x)
\mathbf E_{x,x}^t\left[ {\rm e}^{-\int_0^t V(X_s,\omega)\,{\rm d} s} \mathbf{1}_{\{\tau_{\Delta}>t\}}\right]\,{\rm d} m(x),
\end{equation}
so that $L_M^D=L^{\mathcal G_M}.$

In this section, we will work under the additional assumption that the Laplace exponent $\phi$ of the subordinator $S$ is a complete Bernstein function satisfying the following condition:
\begin{itemize}
\item[\bf(L1)] there exist constants
 $c_{3.1}>0,$ $\beta\in(0,d_w]$ and $s_0>0$ such that for $s \in (0,s_0]$ one has $\phi(s)\leq~c_{3.1} s^{\beta/d_w}.$
\end{itemize}
Under {\bf (L1)}, the assumption {\bf(S2)} is automatically satisfied (it follows e.g. from Lemma \ref{lem:tails}).

\smallskip

We have the following lower bound.

\begin{theo}\label{th:lower}  Let $X$ be a subordinate Brownian motion on $\cG$ via a complete subordinator $S$ with Laplace exponent $\phi$ such that {\bf (S1)} and {\bf(L1)} hold and let $V$ be a Poissonian potential with  profile $W$ satisfying \textbf{(W1)}--\textbf{(W3)}. Suppose that there exist constants $\theta >0,  $  $K \in [0,\infty)$ such that
 \begin{equation}\label{eq:cond-W}
 \limsup_{d(x,y)\to\infty} W(x,y) d(x,y)^{d+\theta} =K.\end{equation}
Then there exist positive constants $C_1,C_2$ and $t_0>0$ such that for $t>t_0$  one has
\begin{equation}\label{eq:lower-0}
L(t)\geq \exp\left\{-C_1 t^{d/(d+\beta)}\nu^{\beta/(d+\beta)} - C_2t^{d/(d+\theta)}\nu\right\}.
\end{equation}
In particular, \\
\p (i) when $\beta<\theta$ then
\[\liminf_{t\to\infty}\frac{\log L(t)}{t^{d/(d+\beta)}}\geq -C_1 \nu^{\beta/(d+\beta)},\] \\
\p (ii) when $\beta=\theta$ then
\[\liminf_{t\to\infty}\frac{\log L(t)}{t^{d/(d+\beta)}}\geq -C_1 \nu^{\beta/(d+\beta)}-C_2\nu,\]
\p (iii) when $\beta >\theta$ then
\[\liminf_{t\to\infty}\frac{\log L(t)}{t^{d/(d+\theta)}}\geq -C_2 \nu.\]
\end{theo}

\noindent{\bf Proof.}
 For given $a>0,$ let the 'short range' and the `long range' profiles
be given by
\begin{eqnarray}\label{eq:wu-a}
W_a(x,y)= W(x,y)\mathbf 1_{\{d(x,y)\leq a\}},  &\mbox{and} & W^a(x,y)= W(x,y)\mathbf 1_{\{d(x,y)>a\}},
\end{eqnarray}
then let  $V_a,$ $V^a$ be the Poissonian potentials based on $W_a$, $W^a,$ accordingly.
Moreover, for $a>0$  let
\[S_W(a):= \sup_{x\in \mathcal G}\int_\mathcal G W^a(x,y)\,{\rm d}m(y)=\sup_{x\in\mathcal G}\int_{d(x.y)>a}W(x,y)\,{\rm d}m(y).\]
We start with the following estimate, which is valid for $W$ satisfying {\bf (W1)} --{\bf (W3)} and $S-$a complete subordinator satisfying {\bf (S1)} and {\bf (S2)}. We prove
that
for any $t>0,$ $a>0,$ $M\in \mathbb Z_+$ one has
\begin{equation}\label{eq:lower-1}
L(t)\geq \exp\left\{ -t\phi\left(\frac{1}{2^{Md_w}}\lambda_1^{BM}(\mathcal G_0)\right) -\nu\,t\,S_W(a)- \nu(2^{Md} +9a^d)\right\}.
\end{equation}
To prove (\ref{eq:lower-1}),
consider the expressions $L^D_{M+k}(t,\omega) = L^{\mathcal G_{M+k}}(t,\omega),$ $k \in \mathbb Z_+$.
Clearly, by (\ref{eq:L-lim}), we have
\[L(t)=\lim_{k\to\infty} \mathbb {E_Q}[L^D_{M+k}(t,\omega)].\]
For given $k\geq 1,$ the set $\mathcal G_{M+k}$ consists of $3^k$ gasket triangles of size $2^{M}$ each, with pairwise disjoint interiors. Denote them $\Delta_1,...,\Delta_{3^k}.$ Because of the inclusions $\Delta_i\subset \mathcal G_{M+k},$ one has
\begin{eqnarray*}
\mathbb{E_Q}[L^D_{M+k}(t,\omega)] & = & \frac{1}{3^{M+k}}\int_{\mathcal G_{M+k}} p(t,x,x) \mathbb {E_Q}\mathbf E^t_{x,x}\left[{\rm e}^{-\int_0^t V(X_s,\omega){\rm d} s}\mathbf 1_{\{\tau_{\mathcal G_{M+k}}>t\}}\right]\,{\rm d}m(x)\\[2mm]
&=&\frac{1}{3^{M+k}}\sum_{i=1}^{3^{k}}\int_{\Delta_i}
p(t,x,x) \mathbb {E_Q}\mathbf E^t_{x,x}\left[{\rm e}^{-\int_0^t V(X_s,\omega){\rm d} s}\mathbf 1_{\{\tau_{\mathcal G_{M+k}}>t\}}\right]\,{\rm d}m(x)\\[2mm]
&\geq & \frac{1}{3^{M+k}}\sum_{i=1}^{3^{k}}\int_{\Delta_i}
p(t,x,x) \mathbb {E_Q}\mathbf E^t_{x,x}\left[{\rm e}^{-\int_0^t V(X_s,\omega){\rm d} s}\mathbf 1_{\{\tau_{\Delta_i}>t\}}\right]\,{\rm d}m(x)\\[2mm]
&\geq&\inf_i \mathbb{E_Q}[L^{\Delta_i}(t,\omega)].
\end{eqnarray*}
Pick now any of the $i$'s, say $i_0,$ and let
\(\mathcal M^a_{i_0}=\{\omega:\mbox{no Poisson points fell
into $\Delta_{i_0}^a$}\},\) where $\Delta_{i_0}^a$ denotes the $a-$vicinity of $\Delta_{i_0}$ . In particular,
\begin{equation}\label{eq:lower-3}\mathbb{E_Q}L^{\Delta_{i_0}}(t,\omega)\geq \mathbb{E_Q}\left[L^{\Delta_{i_0}}(t,\omega)\mathbf 1_{\mathcal M^a_{i_0}}\right].\end{equation}
Observe that for every $\omega \in \mathcal M^a_{i_0}$ and for a fixed trajectory of $X_s$ starting at $x\in \Delta_{i_0}$ and not leaving the set $\Delta_{i_0}$ up to time $t$ one has
\begin{equation}
V(X_s,\omega)= \int_{(\Delta_{i_0}^a)'} W(X_s,y)\,d\mu^{\omega}(y)=
\int_{(\Delta_{i_0}^a)'} W^a(X_s,y)\,{\rm d}\mu^\omega(y)=V^a(X_s,\omega).
\end{equation}
For such a trajectory, the random Feynman-Kac functional ${\rm e}^{-\int_0^t V^a(X_s,\omega)\,{\rm d}s}$ and the event $\mathcal M^a_{i_0}$ are $\mathbb Q-$independent. Therefore, on the set $\{\tau_{\Delta_{i_0}}>t\}$ one has:
\begin{eqnarray}\label{eq:pot1}
\mathbb{E_Q}\left[{\rm e}^{-\int_0^t V(X_s,\omega){\rm d} s}\mathbf 1_{\mathcal M^a_{i_0}}\right]=
\mathbb{E_Q}\left[{\rm e}^{-\int_0^t V^a(X_s,\omega)\,{\rm d}s}\mathbf 1_{\mathcal M_{i_0}^a}\right]=
\mathbb{E_Q}\left[{\rm e}^{-\int_0^t V^a(X_s,\omega){\rm d} s}\right]\cdot\mathbb Q[\mathcal M^a_{i_0}].
\end{eqnarray}

Using the definition (\ref{eq:L-delta}) of $L^{\Delta_{i_0}}$, then inserting (\ref{eq:pot1}) inside (\ref{eq:lower-3})
we obtain:
\begin{eqnarray}\label{eq:pot2}
\mathbb{E_Q}[L^{\Delta_{i_0}}(t,\omega)] &\geq &
\frac{1}{m(\Delta_{i_0})} \int_{\Delta_{i_0}} p(t,x,x) \mathbf E^t_{x,x}\left[
\mathbf 1_{\{\tau_{\Delta_{i_0}}>t\}}\cdot\mathbb{E_Q}\left[{\rm e}^{-\int_0^t V^a(X_s,\omega){\rm d}s}\right]\cdot\mathbb Q[\mathcal M^a_{i_0}]\right]\,{\rm d} m(x),
\nonumber\\
\end{eqnarray}
moreover
\begin{eqnarray*}
\mathbb Q[\mathcal M^a_{i_0}] &=&\mathbb Q[\mbox{no Poisson points inside $\Delta_{i_0}^a$}]
= {\rm e}^{-\nu  m(\Delta_{i_0}^a)}\geq {\rm e}^{-\nu(2^{Md}+9a^{d})}.\end{eqnarray*}

The exponential formula (\ref{eq:exp-formula}) applied to the inner integral in (\ref{eq:pot2})
gives
\[\mathbb{E_Q} {\rm e}^{-\int_0^t V^a(X_s,\omega)\,{\rm d}s} = {\rm e}^{-\nu \int_\mathcal G (1-{\rm e}^{-\int_0^t W^a(X_s,y){\rm d}s}){\rm d}m(y)}.\]
From an elementary inequality $e^{-x}\geq 1-x, $ $x\in \mathbb R,$  and the Fubini theorem we obtain that
\[\mathbb{E_Q} {\rm e}^{-\int_0^t V^a(X_s,\omega)\,{\rm d}s} \geq
{\rm e}^{ -\nu  \int_\mathcal G \int_0^t W^a(X_s, y)\,{\rm d}m(y)\,{\rm d}s} \geq {\rm e}^{-\nu t \sup_{x\in\mathcal G} \int_{d(x,y)> a} W(x,y)\,{\rm d}m(y)} = {\rm e}^{-\nu t S_W(a)}.\]
It follows
\begin{eqnarray*}
\mathbb{E_Q}\left[L^{\Delta_{i_0}}(t,\omega)\right] \geq
\left[ \frac{1}{m(\Delta_{i_0})}\int_{\Delta_{i_0}}p(t,x,x)\mathbf P_{x,x}^t\left[\tau_{\Delta_{i_0}}>t\right]\,{\rm d}m (x)\right]\cdot {\rm e}^{-\nu t S_W(a)}\cdot {\rm e}^{-\nu(2^{Md}+9a^d)}
\end{eqnarray*}

The first multiplier in the expression above
is the averaged trace of the operator $P_t^{\Delta_{i_0}},$ and as such is not bigger than ${\rm e}^{-t\lambda_1(\Delta_{i_0})}.$
From \cite[Theorem 3.4]{bib:CS} we have $\lambda_1(\Delta_{i_0})\leq \phi (\lambda^{BM}_1(\Delta_{i_0})),$ where $\lambda^{BM}_1(U)$ denotes the principal Dirichlet eigenvalue of the Brownian motion on $U.$  As the Brownian motions on $\Delta_{i_0}$
and on $\mathcal G_{M}$ are indistinguishable up to respective exit times of $\Delta_{i_0},$ $\mathcal G_{M},$  one has
$\lambda_1^{BM}(\Delta_{i_0})=\lambda_1^{BM}(\mathcal G_{M}),$ and from the Brownian scaling we have $\lambda_1^{BM}(\mathcal G_{M})=\lambda^{BM}_1(2^{M}\mathcal G_0)= \frac{1}{2^{Md_w}}\lambda^{BM}_1(\mathcal G_0).$
Collecting all the estimates above
we obtain the statement
\[
\mathbb {E_Q}[L^D_{M+k}(t,\omega)]\geq \exp\left\{-t\phi\left(\frac{1}{2^{Md_w}}\lambda_1^{BM}(\mathcal G_0) \right) -\nu\,t\,S_W(a)- \nu(2^{Md} +9a^d)\right\},
\]
which after letting $k$ go to infinity  gives (\ref{eq:lower-1}).

Having proven  (\ref{eq:lower-1}), we will now use the remaining assumptions. From {\bf(L1)} we get
\begin{equation}\label{eq:lower-5}
\phi\left(\frac{1}{2^{Md_w}}\lambda_1^{BM}(\mathcal G_0) \right)\leq \frac{c_{3.1}}{2^{M\beta}}\cdot \left(\lambda_1^{BM}(\cG_0)\right)^{\beta/d_w}=:\frac{c_{3.2}}{2^{M\beta}},
\end{equation}
for $M$ large enough.

The condition (\ref{eq:cond-W}) combined with (\ref{eq:meas-int}) permit us to write
\begin{equation}
\int_{d(x,y)>a} W(x,y)\,{\rm d}m(y) \leq  \int_{d(x,y)>a}\frac{(K+o(1))}{{d(x,y)}^{d+\theta}} {\rm d}m (y)
\leq c_{2.2}(K+o(1)) \frac{1}{a^\theta}, \nonumber
\end{equation}
and consequently
\begin{equation}\label{eq:es-a}
S_W(a)\leq c_{2.2}(K+o(1)) \frac{1}{a^\theta}, \quad \mbox{as} \quad a \to \infty.
\end{equation}

Next, for given $t>0,$  choose $M=M(t)$ to be the unique integer satisfying
\begin{equation}\label{eq:M-zero}
2^{M}\leq \left(\frac{t}{\nu}\right)^{\frac{1}{d+\beta }}< 2^{M+1}.
\end{equation}
  Inserting (\ref{eq:lower-5}), (\ref{eq:es-a}), and (\ref{eq:M-zero}) into
(\ref{eq:lower-1}), and using $a= t^{1/(d+\theta)},$ after some elementary algebra we obtain
\[L(t)\geq \exp\left\{-\left({c_{3.2}}+1\right) t^{d/(d+\beta)}\nu^{\beta/(d+\beta)} - (Kc_{2.2}+9+o(1)) t^{d/(d+\theta)}\nu\right\}, \quad \mbox{as} \quad t \to \infty.
\]
To get (\ref{eq:lower-0}), we just set
$C_1= \left(c_{3.2}+1\right),$ $C_2= (Kc_{2.2} +9).$ Statements (i)--(iii) are straightforward consequences of
(\ref{eq:lower-0}). \hfill$\Box$

\subsection{Lower bounds for the Feynman-Kac functional}
The methods we use are also suitable for obtaining bounds on the averaged Feynman-Kac functional, i.e. $\mathbb{E_Q}\mathbf E_x\left[{\rm e}^{-\int_0^t V(X_s,\omega)\,{\rm d}s}\right].$
In this case, the assumptions concerning the process and the profile $W$ can be somewhat relaxed.

\begin{theo}\label{th:lower-funct-1} Let $X$ be a subordinate Brownian motion on $\cG$ via a complete subordinator $S$ with Laplace exponent $\phi$ such that {\bf (S1)} and {\bf(L1)} hold  and let the  potential profile  $W$ fulfil  {\bf (W1)}  and (\ref{eq:cond-W}).
Then, with constants $C_1,C_2>0$ from Theorem \ref{th:lower}, for any $x\in\mathcal G$ one has:\\
\p (i) when $\beta<\theta$ then
\[\liminf_{t\to\infty}\frac{\log \mathbb{E_Q}\mathbf E_x[{\rm e}^{-\int_0^t V(X_s,\omega)\,{\rm d}s}]}{t^{d/(d+\beta)}}\geq -C_1 \nu^{\beta/(d+\beta)},\] \\
\p (ii) when $\beta=\theta$ then
\[\liminf_{t\to\infty}\frac{\log  \mathbb{E_Q}\mathbf E_x[{\rm e}^{-\int_0^t V(X_s,\omega)\,{\rm d}s}]}{t^{d/(d+\beta)}}\geq -C_1 \nu^{\beta/(d+\beta)}-C_2\nu,\]
\p (iii) when $\beta >\theta$ then
\[\liminf_{t\to\infty}\frac{\log  \mathbb{E_Q}\mathbf E_x[{\rm e}^{-\int_0^t V(X_s,\omega)\,{\rm d}s}]}{t^{d/(d+\theta)}}\geq -C_2 \nu.\]
\end{theo}

\
\noindent {\bf Proof.}
The proofs of these statements are very much alike those from Theorem  \ref{th:lower}. Fix $M\in\mathbb Z_+,$  $a>0,$ $t>0,$ $x\in\mathcal G$, for the moment assuming only that $M$ is so large that $x\in\mbox{Int}\,\mathcal G_M.$ Let $\mathcal M_M^{{a}}$ be the event `no Poisson points fell into $\mathcal G_M^a$'. Using the reasoning that led to (\ref{eq:pot1}),
we get that
\begin{eqnarray*}
\mathbb{E_Q}\mathbf E_x\left[{\rm e}^{-\int_0^t V(X_s,\omega){\rm d} s}\right] &\geq & \mathbb{E_Q}\mathbf E_x\left[{\rm e}^{-\int_0^t V(X_s,\omega){\rm d}s}\mathbf 1_{\{\tau_{\mathcal G_M}>t\}}\mathbf 1_{\mathcal M_M^{a}}\right]\\
&\geq & \mathbf E_x\left[\mathbf 1_{\{\tau_{\mathcal G_M}>t\}}\mathbb {E_Q}\left[{\rm e}^{-\int_0^t V^a(X_s,\omega){\rm d}s}\right]\right]\mathbb Q[\mathcal M_M^{a}]\\
&\geq & \mathbf P_x[\tau_{\mathcal G_M}>t]\cdot {\rm e}^{-\nu t S_W(a) -\nu (2^{Md} +9a^d)}.
\end{eqnarray*}

The probability space $(\Omega, \cF, \pr_x)$ for $X$ can be realized as a product of two probability spaces on which $Z$ and $S$ are defined. In particular, $\pr_x = \pr^X_x := \pr^Z_x \otimes \pr$. Therefore, it is easy to observe that for every  $M \in \mathbb Z_+,$ $x \in \cG_M$,  and $t>0$, by Fubini we have
$$
\pr^X_{x} [ \tau^X_{\cG_M}>t] \geq \pr^Z_{x} \otimes \pr[ \tau^Z_{\cG_M}>S_t] = \int_0^{\infty} \pr^Z_{x}[\tau^Z_{\cG_M}>u ]\ \eta_t({\rm d}u).
$$
By the scaling of $Z$, we have $\pr^Z_{x}[ \tau^Z_{\cG_M}>u] =\pr^Z_{(x/2^M)}[\tau^Z_{\cG_0}> 2^{-Md_w} u  ]$. Moreover, one has
$$
c \pr^Z_x[\tau^Z_{\cG_0}> u ] \geq \pr^Z_x[\tau^Z_{\cG_0}> u; \varphi_1^{0,Z}(Z_t)] = e^{-u \lambda_1^{BM}(\mathcal G_0) } \varphi_1^{0,Z}(x)  , \quad x \in \cG_0, \ \ u>0,
$$
where $\lambda_1^{BM}(\mathcal G_0)$ is the principal eigenvalue of the Brownian motion on $\cG_0$ with killing on $\partial \cG_0,$ $\varphi_1^{0,Z}$ is the corresponding normalized eigenfunction, and $c = \|\varphi_1^{0,Z}\|_{\infty} < \infty$ is independent of $M$ and $u$. The transition density of the process $Z$ killed on exiting $\mathcal G_0$ is positive for all $u>0$, $x,y\in \mbox{Int}\,\mathcal G_0,$  thus
from general theory its ground state is continuous and can be chosen to be strictly positive on $\overline{B}(0,\frac{1}{2})$ (note  that $0\notin \partial\mathcal G_0\subset\mathcal G$). Denote $c^{(1)}=\inf_{y\in B(0,1/2)} \varphi_1^{0,Z}(y) >0$ and $c^{(2)}=c^{-1} \cdot c^{(1)}$. Collecting all the above estimates, we get
$$
\pr^X_{x} [\tau^X_{\cG_M}>t] \geq \int_0^{\infty} \pr^Z_{(x/2^M)}[ \tau^Z_{\cG_0}>2^{-Md_w} u ] \ \eta_t({\rm d}u) \geq c^{(2)} \int_0^{\infty} e^{-2^{-Md_w} u \lambda_1^{BM}(\mathcal G_0) } \ \eta_t({\rm d}u),
$$
whenever $x \in \cG_{M}/2$, i.e., $x/2^M \in B(0,1/2)$. Observe the last integral is the Laplace transform of $\eta_t$. Therefore,
$$
\pr^X_{x} [ \tau^X_{\cG_M}>t] \geq c^{(2)} {\rm e}^{-t \phi(2^{-Md_w} \lambda_1^{BM}(\mathcal G_0))}.
$$
 Using now the condition {\bf(L1)}, we finally obtain
\[\mathbf P_x^X[\tau^X_{\mathcal G_M}>t]\geq c^{(2)} {\rm e}^{-c_{3.1} t\cdot 2^{-M\beta d_w}}. \]
All these arguments were true for any $M,$ as long as $x\in\mathcal G_M/2.$ At this point we declare the specific value of $M$: we take it equal to $M(t)$ given by (\ref{eq:M-zero}). Moreover, choose again $a=t^{\frac{1}{d+\theta}}.$
For this choice of $M,a,$ identically as before we obtain, as long as $x\in \mathcal G_{M(t)}/2$ and $t \to \infty$,
\[\mathbb{E_Q}\mathbf E_x\left[{\rm e}^{-\int_0^t V(X_s,\omega)\,{\rm d}s}\right]\geq c^{(2)} \exp\left\{-\left(c_{3.2}+1\right) t^{d/(d+\beta)}\nu^{\beta/(d+\beta)} - (Kc_{2.2}+9+o(1)) t^{d/(d+\theta)}\nu\right\}.
\]
From this inequality all the statements follow as before.

\hfill$\Box$

\section{Upper bounds}

\subsection{Upper bound for the long range interaction}

We first derive the upper bound which depends only on the long range behaviour of the potential. It does not require any additional assumptions on the subordinator $S$. The following result is useful for profile functions $W$ with slow decay at infinity.

\begin{prop} \label{prop:long_range_2} Let $X$ be a subordinate Brownian motion on $\cG$ via a subordinator $S$ satisfying \textbf{(S1)}-\textbf{(S2)} and let $V$ be a Poissonian potential with  profile $W$ such that the assumptions \textbf{(W1)}-\textbf{(W3)} hold.
Then for every $t \geq 1$ and $a>0$ we have
\begin{equation} \label{eq:bound1}
L(t) \leq c_{2.9} e^{-\nu R_W(a,t)},
\end{equation}
where $R_W(a,t):= \inf_{x \in \cG} \int_{d(x,y)>a} \left(1-e^{-t W(x,y)}\right) m({\rm d} y)$, and $c_{2.9}$ is the constant from (\ref{eq:sup-of-g}). In particular, if for a number $\theta>0$ there is $K \in [0,\infty)$ such that
\begin{equation} \label{eq:cond-W_up}
\liminf_{d(x,y) \to \infty} W(x,y) d(x,y)^{d+\theta} = K,
\end{equation}
then we have
$$
\limsup_{t \to \infty} \frac{\log L(t)}{t^{\frac{d}{d+\theta}}} \leq - \frac{K}{c_{2.2} e^K} \nu.
$$
\end{prop}

\noindent{\bf Proof.}
Since for every $t>0$ we have
$L(t) = \lim_{M \to \infty} \mathbb{E}_{\qpr} L_M^D(t,\omega)$, it is enough to show that for every $t \geq 1$, $a>0$ we have
$$
\mathbb{E}_{\qpr} L_M^D(t,\omega) \leq c_{2.9} e^{-\nu R_W(a,t)}.
$$
Recall that for every $t \geq1$ and $x ,y \in \cG$ we have $p(t,x,y) \leq c_{2.9} $ (see (\ref{eq:sup-of-g})). By this bound and the exponential formula (\ref{eq:exp-formula}), we thus get
\begin{eqnarray*}
\mathbb{E}_{\qpr} L_M^D(t,\omega) & \leq & \frac{c_{2.9}}{m(\cG_M)}\int_{\cG_M}
\mathbf E_{x,x}^t \mathbb{E}_{\qpr} \left[{\rm e}^{-\int_0^t V(X_s,\omega)\,{\rm d} s}\right]\,m({\rm d}x) \\ \nonumber
& = & \frac{c_{2.9}}{m(\cG_M)}\int_{\cG_M}
\mathbf E_{x,x}^t\left[  {\rm e}^{- \nu \int_{\cG} \left(1-e^{-\int_0^t W(X_s,y)}{\rm d} s\right) m({\rm d} y)}\right]\,m({\rm d} x).\\ \nonumber
\end{eqnarray*}
The integral over $\cG$ in the exponent can be written as
\begin{eqnarray*}
\int_{\cG} \left(1-e^{-\int_0^t W(X_s,y)}{\rm d} s\right) m({\rm d} y) &= &
  \int_{\cG}F \left(\int_0^t t W(X_s,y)\frac{{\rm d}s}{t}\right)\,m({\rm d} y),
\end{eqnarray*}
where the function $F$ is given by $F(t)=1-e^{-t}.$ It is a
concave function, therefore from Jensen's inequality for concave
functions and Fubini's theorem we get
\begin{eqnarray*}
\int_{\cG} \left(1-e^{-\int_0^t W(X_s,y)}{\rm d} s\right) m({\rm d} y) & = & \int_{\cG}  F \left(\int_0^t t W(X_s,y)\frac{{\rm d}s}{t} \right) m({\rm d} y) \\
&\geq &
\int_{\cG} \int_0^t F \left(t W(X_s,y)\right) \frac{{\rm d}s}{t}m({\rm d} y) \\
&=& \int_0^t \int_{\cG}  \left(1-e^{-t W(X_s,y)} \right)  m({\rm d} y) \frac{{\rm d}s}{t} \\
& \geq & \inf_{x \in \cG} \int_{\cG}  \left(1-e^{-t W(x,y)}\right)  m({\rm d} y).
\end{eqnarray*}
In particular, for all $t \geq 1$, $a>0$ and $M \in \mathbb Z_+$, we obtain
\begin{eqnarray*}
\mathbb{E}_{\qpr}[ L_M^D(t,\omega) ]\leq c_{2.9} e^{- \nu \inf_{x \in \cG} \int_{d(x,y)>a}  \left(1-e^{-t W(x,y)}\right)  m({\rm d} y)}=c_{2.9} {\rm e}^{-\nu R_W(a,t)},
\end{eqnarray*}
which completes the proof of (\ref{eq:bound1}).

To show the second assertion, first note that by the standard inequality $1-e^{-s} \geq e^{-s} s$, $s \geq 0$, and (\ref{eq:cond-W_up}), for $d(x,y) > a$ and $t > 0$, we have
$$
1-e^{-tW(x,y)} \geq 1 - e^{-t (K+o(1)) \, d(x,y)^{-d-\theta}} \geq t (K+o(1)) \, d(x,y)^{-d-\theta} e^{-t (K+o(1)) \, d(x,y)^{-d-\theta}} \quad \mbox{as} \quad a \to \infty.
$$
By taking $a = t^{\frac{1}{d+\theta}}$ with $t \to \infty$, we thus get, using (\ref{eq:meas-int}):
\begin{eqnarray*}
\int_{d(x,y)>t^{\frac{1}{d+\theta}}} \left(1-e^{-tW(x,y)}\right) m({\rm d} y) & \geq & t (K+o(1)) e^{-(K+o(1)) } \int_{d(x,y)>t^{\frac{1}{d+\theta}}} d(x,y)^{-d-\theta} m({\rm d} y) \\
& \geq & t (K+o(1)) e^{-(K+o(1)) } c^{-1}_{2.2} t^{\frac{-\theta}{d+\theta}} \\
& = & c^{-1}_{2.2}t^{\frac{d}{d+\theta}} (K+o(1)) e^{-(K+o(1)) }.
\end{eqnarray*}
By (\ref{eq:bound1}), we conclude that
$$
\limsup_{t \to \infty} \frac{\log L(t)}{t^{{d}/{d+\theta}}} \leq - \frac{K}{c_{2.2}e^K} \nu.
$$
The proof is complete.

\hfill$\Box$

Similar bounds hold true for the averaged Feynman-Kac functional.

\begin{prop} \label{prop:long_range_funk_2} Let $X$ be a subordinate Brownian motion via the subordinator $S$  satisfying {\bf (S1)} -- {\bf(S2)} and let $V$ be a Poissonian potential with  profile $W$ satisfying \textbf{(W1)}.
Then for every $t \geq 1$ and $a>0$ we have
\begin{equation} \label{eq:bound2}
\sup_{x \in \cG} \mathbb{E}_{\qpr}  \ex_x\left[e^{-\int_0^t V(X_s,\omega) ds}\right] \leq  e^{-\nu R_W(a,t)}.
\end{equation}
In particular, if (\ref{eq:cond-W_up}) holds with some $\theta>0$ and $K \in [0,\infty)$,
then for every $x\in\mathcal G$ we have
$$
\limsup_{t \to \infty} \frac{\log \left( \mathbb{E}_{\qpr} \ex_x\left[e^{-\int_0^t V(X_s,\omega) ds}\right]\right)}{t^{{d}/{d+\theta}}} \leq - \frac{K}{c_{2.2} e^K} \nu.
$$
\end{prop}

\noindent{\bf Proof.} For every $x \in \cG$, we have
$$
\mathbb{E}_{\qpr}  \ex_x\left[e^{-\int_0^t V(X_s) ds}\right] = \ex_x\left[{\rm e}^{- \nu \int_{\cG} \left(1-e^{-\int_0^t W(X_s,y)}{\rm d} s\right) m({\rm d} y)}\right] \leq e^{-\nu R_W(a,t)}, \quad a>0,
$$
and the second assertion follows exactly by the same argument as in Proposition \ref{prop:long_range_2}.
\hfill$\Box$

\subsection{Upper bound for the short range interaction}

Recall that we assume the Laplace exponent $\phi$ to be a complete Bernstein function of the form
\begin{equation} \label{eq:def_phi}
\phi(\lambda) = b \lambda + \psi(\lambda) \quad \mbox{with} \quad \psi(\lambda)= \int_0^{\infty} \left(1-e^{-\lambda u} \right) \rho(u) {\rm d}u, \quad \lambda \geq 0.
\end{equation}
In this subsection we need stronger  assumptions on the Laplace exponent $\phi$:
\begin{itemize}
\item[\textbf{(U1)}] $b > 0$ and $\psi \equiv 0$ (equivalently, $\nu \equiv 0;$ no jumps)

or

\item[\textbf{(U2)}] $b > 0$ and $\psi \neq 0$ satisfies the following weak scaling conditions: there are $\alpha_1,\alpha_2, \beta, \delta \in (0,d_w)$, $a_1, a_2 \in (0,1]$, $a_3, a_4 \in [1,\infty)$ and $r_0 >0$ such that
\begin{eqnarray} \label{eq:small_scaling}
&&a_1 \lambda^{\alpha_1/d_w} \psi(r) \leq \psi(\lambda r) \leq a_3 \lambda^{\beta/d_w} \psi(r), \quad \lambda \in (0,1], \quad r \in (0, r_0]
\\[2mm]
\label{eq:large_scaling} \!\!\!\!\!\!\! \mbox{and}\;\;\;\;\;\;\;\;\;\;\;\;&&
a_2 \lambda^{\alpha_2/d_w} \psi(r) \leq \psi(\lambda r) \leq a_4 \lambda^{\delta/d_w} \psi(r), \quad \lambda \geq 1, \quad r \geq r_0
\end{eqnarray}

or

\item[\textbf{(U3)}] $b =  0$ and $\psi \neq 0$ satisfies (\ref{eq:small_scaling}) and (\ref{eq:large_scaling}) with
    $\alpha_1=\alpha_2.$
\end{itemize}
Note that under the assumption \textbf{(U1)} the subordinator $S$ is a pure drift, while the left hand sides of (\ref{eq:small_scaling}) and (\ref{eq:large_scaling})  imply the lower bounds
\begin{eqnarray} \label{eq:small_low_scaling}
\phi(\lambda) \geq \psi(\lambda) \geq \overline{a}_1 \ \lambda^{\alpha_1/d_w}, \quad \lambda \in(0,1],
\\
\label{eq:large_low_scaling}
\phi(\lambda) \geq \psi(\lambda) \geq \overline{a}_2 \ \lambda^{\alpha_2/d_w}, \quad \lambda \in[1,\infty)
\end{eqnarray}
(we have set $\overline{a}_i=a_i\psi(r_0)r_0^{-\alpha_i/d_w}).$

Moreover, one can directly check that if \textbf{(U1)}, \textbf{(U2)},
 or \textbf{(U3)} is satisfied, then both assumptions \textbf{(S1)} and \textbf{(S2)} hold (see \cite[Remark 2.1 (2) and Lemma 2.2]{bib:KaPP}). Assumption {\bf (L1)} is in this case true as well.
Examples of subordinators with Laplace exponents satisfying  {\bf (U1)} -- {\bf (U3)} and the corresponding subordinate Brownian motions will be discussed in  Section \ref{sec:examples}.

We will need the following estimates, which are  consequences of \textbf{(U2)} or \textbf{(U3)}.
\begin{lem} \label{lem:comp_levy} Let $S$ be a complete subordinator with  Laplace exponent $\phi$ given by (\ref{eq:def_phi}).
Under the condition \textbf{(U2)} or \textbf{(U3)} the following estimates hold.
\begin{itemize}
\item[(a)] There exists a constant $c_{4.1}=c_{4.1}(\phi) \in (0,1]$ such that
$$
\rho(s) \geq c_{4.1} \, s^{-1} \cdot \left\{
  \begin{array}{ccc}
    s^{-\alpha_1/d_w} & \mbox{  if  } & s \geq 1,\\
    s^{-\alpha_2/d_w} & \mbox{  if  } & s \in(0,1].
  \end{array}\right.
$$
\item[(b)] There exists a constant $c_{4.2}=c_{4.2}(\phi) > 0$ such that
$$
\int_0^{\infty} u^{-d_s/2} \eta_t({\rm d}u) \leq c_{4.2} \, \left(t^{-d/\alpha_1}+ t^{-d/\alpha_2}\right), \quad t > 0.
$$
\end{itemize}
\end{lem}

\noindent{\bf Proof.}
We first prove (a). By \cite[Proposition 2.5]{bib:KSV} the conditions (\ref{eq:small_scaling}) and (\ref{eq:large_scaling}) imply that there is a constant $c^{(1)}>0$ such that
$$
\rho(s) \geq c^{(1)} s^{-1} \psi(s^{-1}), \quad s >0.
$$
This, together with (\ref{eq:small_low_scaling}) and (\ref{eq:large_low_scaling}), imply the claimed inequalities in (a).

Consider now (b). We have
$$
e^{-t\phi(\lambda^{2/d_s})} = \int_0^{\infty} e^{-\lambda^{2/d_s} u} \eta_t({\rm d} u), \quad t>0, \ \lambda>0.
$$
By integrating in this equality with respect to $\lambda$ over $(0,\infty)$ and by Fubini, we get
$$
\int_0^{\infty} e^{-t\phi(\lambda^{2/d_s})} {\rm d} \lambda = \int_0^{\infty} \left(\int_0^{\infty}e^{-(u^{d_s/2}\lambda )^{2/d_s}} {\rm d} \lambda \right) \eta_t({\rm d} u), \quad t>0.
$$
Now, the substitution $\vartheta = u^{d_s/2}\lambda$ in the internal integral on the right hand side gives
$$
\int_0^{\infty} e^{-t\phi(\lambda^{2/d_s})} {\rm d} \lambda = \int_0^{\infty}e^{-\vartheta^{2/d_s}} {\rm d} \vartheta \cdot \int_0^{\infty} u^{-d_s/2} \eta_t({\rm d} u), \quad t>0,
$$
that is,
$$
\int_0^{\infty} u^{-d_s/2} \eta_t({\rm d} u) = \frac{1}{c^{(2)}} \int_0^{\infty} e^{-t\phi(\lambda^{2/d_s})} {\rm d} \lambda, \quad t>0,
$$
with $(0,\infty) \ni c^{(2)} := \int_0^{\infty}e^{-\vartheta^{2/d_s}} {\rm d} \vartheta$. It is enough to estimate the integral on the right hand side. Recalling that $d_s=2d/d_w$ and applying the bounds (\ref{eq:small_low_scaling}), (\ref{eq:large_low_scaling}),
$$
\int_0^{\infty} e^{-t\phi(\lambda^{2/d_s})} {\rm d} \lambda \leq \left( \int_0^{1} e^{-\overline {a}_1 t \lambda^{\alpha_1/d}} {\rm d} \lambda + \int_1^\infty e^{-\overline {a}_2 t \lambda^{\alpha_2/d}} {\rm d} \lambda\right), \quad t>0.
$$
Finally, by substitution $\vartheta = t^{d/\alpha_i}\lambda$ in respective integrals, we conclude that
$$
\int_0^{\infty} e^{-t\phi(\lambda^{2/d_s})} {\rm d} \lambda \leq c^{(3)} \left( t^{-d/\alpha_1}+t^{-d/\alpha_2}\right), \quad t>0.
$$
We set $c_{4.2}=  c^{(3)}/c^{(2)}.$ The proof is complete.
\hfill$\Box$

\bigskip
\noindent
Observe that by Lemma \ref{lem:comp_levy} (b) we have  $c_{2.7}(t) \leq c_{4.2}\left( t^{-d/\alpha_1}+t^{-d/\alpha_2}\right)$, where $c_{2.7}(t)$ comes from \textbf{(S1)}.

\subsubsection{Reflected subordinate Brownian motions and their Schr\"odinger perturbations}

Our results in this section strongly rely on the so-called reflected subordinate Brownian motions introduced recently in \cite{bib:KaPP}. Therefore first we need to make a necessary preparation. For more detail discussion and justification of all  properties of reflected processes listed below we refer the reader to \cite[Subsection 2.2.3]{bib:KaPP} and references therein.

Let $M \in \mathbb Z_+$ and let $Z^M$ be the reflected Brownian motion in $\cG_M$ introduced in \cite{bib:KPP-PTRF}, i.e. a Feller diffusion with strictly positive transition densities with respect to $m,$ which are given by the formula
\[g^M(t,x,y)=\left\{\begin{array}{ll}
\sum_{y'\in \pi_M^{-1}(y)} g(t,x,y'), & \mbox{when } x,y\in\mathcal G_M,\, y\notin \mathcal V_M\setminus\{0\},\\[2mm]
2\sum_{y'\in \pi_M^{-1}(y)} g(t,x,y'), & \mbox{when } x\in\mathcal G_M,\, y\in \mathcal V_M \backslash \left\{0\right\},\end{array}
\right.\]
where $\pi_M$ is the projection described in Subsection \ref{subsec:gasket}. The function $g^M(t,x,y)$ is jointly continuous in $(t,x,y)$ and symmetric in its space variables. It follows from scaling properties of $g$ and properties of the projections $\pi_M$ that
\begin{equation} \label{eq:scal_gM}
g^M(t,x,y) = 2^{-Md} g^0(2^{-Md_w}t,2^{-M}x, 2^{-M}y), \quad x, y \in \cG_M, \ t >0, \ M \in \mathbb Z_+.
\end{equation}
The transition semigroup of the processes $Z^M$ and the corresponding Dirichlet forms will be denoted by $(G^M_t)_{t \geq 0}$ and $(\cE_{(d_w)}^{M}, \cF_{(d_w)}^{M})$, respectively. Recall that
\begin{eqnarray} \label{eq:dir_form_ZM}
\cE_{(d_w)}^{M}(u,u) & := &\lim_{t \to 0^{+}} \left(\frac{u-G^M_tu}{t}, u\right)_{L^2(\cG_M,m)} \nonumber\\  & = & \lim_{t \to 0^{+}} \frac{1}{2t} \int_{\cG_M \times \cG_M} (u(x)-u(y))^2 \, g^M(t,x,y) m({\rm d}x) m({\rm d}y),
\end{eqnarray}
for all functions $u \in \cF_{(d_w)}^{M}$, i.e. for those functions for which the limit in (\ref{eq:dir_form_ZM}) is finite. One can directly check that by (\ref{eq:scal_gM}) we have
$$
\cE_{(d_w)}^{M}(u,u) = 2^{-Md_w} \cE_{(d_w)}^{0}(u_M,u_M) \quad \mbox{with} \quad u_M(\cdot)=2^{Md/2} u(2^M \cdot), \ M \in \mathbb Z_+.
$$

The symmetric Markov process $X^M = (X_t^M, \pr_x^M)_{t \geq 0, \, x \in \cG_M}$ given by $X^M_t := Z^M_{S_t}$ is called the \emph{reflected subordinate Brownian motion} via the subordinator $S$ in $\cG_M$. Throughout this section we always assume that the subordinator $S$ meets one of the assumptions \textbf{(U1)}, \textbf{(U2)} or \textbf{(U3)}, which means that also the both regularity conditions \textbf{(S1)}--\textbf{(S2)} hold. Processes $Z^M$ and $S$ are always assummed to be stochastically independent and, therefore, the subordination formula
\begin{equation}\label{eq:proj-sub}
p^M(t,x,y)=\int_0^\infty g^M(u,x,y)\eta_t({\rm d}u),\;\;\;x,y\in \mathcal G_M, \;t>0,
\end{equation}
defines the transition densities of the process $X^M$. Kernels $p^M$ inherit the symmetry from $g^M$ and have the same continuity properties as $p$, given by (\ref{eq:subord}). Moreover, when the assumption \textbf{(U1)} holds, we simply have $X^M_t = Z^M_{bt}$ and $p^M(t,x,y)=g^M(bt,x,y)$ for all $t>0$ and $x,y \in \cG_M$, while under \textbf{(U2)}  or \textbf{(U3)}, $X^M$ is a jump process with density $p^M$ satisfying the following upper bound (cf. \cite[formula (2.13)]{bib:KaPP}).

\begin{lem} \label{lem:pMbound} Under the assumption \textbf{(U2)} or \textbf{(U3)}
there is a constant $c_{4.3}=c_{4.3}(\phi) >0$ such that
\begin{equation}\label{eq:pMbound1}
p^M(t,x,y)\leq c_{4.3} \left((t\wedge 2^{M\beta})^{-d/\alpha_1}+
(t\wedge 2^{M\beta})^{-d/\alpha_2}+(t\wedge 2^{M\beta})^{-d/\beta}\right),\quad t>0, \ x, y \in \cG_M, \ M \in \mathbb Z_+.
\end{equation}
\end{lem}

\noindent{\bf Proof.} By lemma \cite[Lemma 2.5, ineq. (2.13)]{bib:KaPP}, we have
$$
g^M(u,x,y) \leq c \left(u^{-d/d_w} \vee 2^{-Md} \right), \quad u>0, \ x, y \in \cG_M, \ M \in \mathbb Z_+,
$$
with an absolute constant $c>0.$ Therefore, by the subordination formula (\ref{eq:proj-sub}), we get
$$
p^M(t,x,y) \leq c\left(\int_0^{2^{Md_w} } u^{-d_s/2} \eta_t({\rm d}u) + 2^{-Md} \eta_t(2^{M d_w}, \infty) \right), \quad t>0, \ x, y \in \cG_M, \ M \in \mathbb Z_+.
$$
Note that by Lemma \ref{lem:comp_levy} (b) the first member of the sum above is smaller than $c^{(1)} \left(t^{-d/\alpha_1}+t^{-d/\alpha_2}\right)$ for some constant $c^{(1)}>0,$ and all $t>0$ and $M \in \mathbb Z_+$. Furthermore, it immediately follows from Lemma \ref{lem:tails} and the upper bound in (\ref{eq:small_scaling}) of \textbf{(U2)} (or \textbf{(U3)}) that
$$
\eta_t(2^{M d_w}, \infty) \leq (c^{(2)}  t 2^{-M \beta}\wedge 1), \quad t >0, \ M \in \mathbb Z_+.
$$
Collecting both estimates above, we obtain
$$
p^M(t,x,y) \leq c^{(3)} \left(t^{-d/\alpha_1} + t^{-d/\alpha_2} + {2^{-Md}}({t}{2^{-M\beta}}\wedge 1)\right), \quad t>0, \ x, y \in \cG_M, \ M \in \mathbb Z_+.
$$
Furthermore, when $t\leq 2^{M\beta},$ then one has ${2^{-Md}}\left({t}{2^{-M\beta}}\wedge 1\right)=
{t}{2^{-M(\beta+d)}}\leq {t}\cdot{t^{-(\beta+d)/\beta}}
={t^{-d/\beta}},
$
while for $t > 2^{M\beta}$ we have ${t^{-d/\alpha_i}}\leq 2^{-M\beta d/\alpha_i}.$
This results in  the bound (\ref{eq:pMbound1}).
\hfill$\Box$

\medskip
By $\pr^{M,t}_{x,y}$ we denote the bridge measures corresponding to process $X^M$ on $D([0,t],\cG_M)$ (for more details we refer to \cite[p. 11-12]{bib:KaPP}).

The process $X^M$ corresponding to the specific subordinator $S$ with Laplace exponent $\phi(\lambda)=\lambda^{\gamma/d_w}$, $\gamma \in (0,d_w]$, will be singled out below. We will denote it by $X^M_{(\gamma)}$ and, by analogy to the Euclidean case, we call it the $\gamma$\emph{-stable reflected subordinate Brownian motion} in $\cG_M$. Clearly, when $\gamma = d_w$, then we just have $X^M_{(\gamma)} = Z^M$. Note that stable processes reflected in $\cG_0$ were recently considered in \cite{bib:Kow-KPP}.

By $(\cE_{\phi}^M, \cF_{\phi}^M)$ we denote the Dirichlet form corresponding to the reflected process $X^M$ in $\cG_M$ (resp. $(\cE_{(\gamma)}^{M}, \cF_{(\gamma)}^{M})$ for $X^M_{(\gamma)}$). We always have $\cF_{(d_w)}^{M} \subset \cF_{\phi}^{M}$. It is known (see \cite{bib:CS,bib:Okura02}) that when $b>0$ then $\cF_{\phi}^{M} = \cF_{(d_w)}^{M},\,$ and for $u \in \cF_{(d_w)}^{M}$ we have
\begin{eqnarray}
\cE_{\phi}^{M}(u,u) & = & b \cE_{(d_w)}^{M}(u,u) + \int_0^{\infty} (u-G^M_s u, u)_{L^2(\cG_M,m)} \rho(s){\rm d} s \nonumber \\
& = & b \cE_{(d_w)}^{M}(u,u) + \int_{\cG_M \times \cG_M} (u(x)-u(y))^2 J_{\phi}^M(x,y) m({\rm d}x) m({\rm d}y), \nonumber
\end{eqnarray}
where
\begin{equation}\label{eq:j-em}
J_{\phi}^M(x,y) = \frac{1}{2} \int_0^{\infty} g^M(s,x,y) \rho(s){\rm d} s.
\end{equation}
For $b=0$ ($S$ has no drift) we have
$$
\cF_{\phi}^M = \left\{ u \in L^2(\cG_M,m): \int_0^{\infty} (u-G^M_s u, u)_{L^2(\cG_M,m)} \rho(s){\rm d} s < \infty \right\}
$$
and for $u \in \cF_{\phi}^M$
$$
\cE_{\phi}^M(u,u) = \int_{\cG_M \times \cG_M} (u(x)-u(y))^2 J_{\phi}^M(x,y) m({\rm d}x) m({\rm d}y).
$$

In the sequel we will investigate the process $X^M$ perturbed by  potentials $V(x)$, $x \in \cG_M$, such that $V \in \cK^{X^M}$. Recall that the Kato class $\cK^{X^M}$ related to the process $X^M$ consists of functions $V$ satisfying the condition $\lim_{t \to 0^{+}} \sup_{x \in \cG_M} \ex^M_x \left[\int_0^t| V|(X_s^M) {\rm d}s\right] = 0$ (one can check that if $V \in \cK_{\loc}^X$, then $V \textbf{1}_{\cG_M} \in \cK^{X^M}$). The corresponding transition semigroup, which we call the Feynman-Kac semigroup associated to the process $X^M$ and the potential $V$, consists of operators
$$
T_t^{\phi, V, M} f(x) = \ex^M_x \left[e^{-\int_0^t V(X_s^M) ds} f(X^M_t) \right], \quad f \in L^2(\cG_M,m), \ t>0.
$$
(for $X^M_{(\gamma)}$, $\gamma \in (0,d_w]$, we write $T_t^{\gamma, V, M}$). Again,
for every $t>0$, the operators $T_t^{\phi, V, M}$ (resp. $T_t^{\gamma, V, M}$) are of Hilbert-Schmidt type and have purely discrete spectrum of the form $\left\{\exp(-t \lambda_n^M(\phi,V))\right\}_{n \geq 1}$ (resp. $\left\{\exp(-t \lambda_n^M(\gamma,V))\right\}_{n \geq 1}$), such that
$0 \leq \lambda_1^M(\phi,V) < \lambda_2^M(\phi,V) \leq \lambda_3^M(\phi,V) \leq ... \to \infty$. For the verification of the above properties and more details on the Feynman-Kac semigroups of the reflected subordinate Brownian motions we refer the reader to \cite[Subsection 2.3.1]{bib:KaPP}.

We also define the Dirichlet form $(\cE_{\phi, V}^M, \cF_{\phi, V}^M)$ corresponding to the 'reflected' process $X^M$ perturbed by $V$ (resp. $(\cE_{(\gamma), V}^{M}, \cF_{(\gamma), V}^{M}),$ for $X^M_{(\gamma)}$). Since $V \in \cK^{X^M}$, we also have $V \in L^1(\cG_M, m)$ and, by general theory of Dirichlet forms \cite[Section 6]{bib:FOT}, it holds that
$$\cF_{\phi, V}^M = \cF_{\phi}^M \cap L^2(\cG_M, V(x)m({\rm d} x))$$
and for $u \in \cF_{\phi, V}^M$ we have
$$
\cE_{\phi, V}^M(u,u)  = \cE_{\phi}^M(u,u) + \int_{\cG_M} V(x) u^2(x) m({\rm d}x).
$$
As above, for $M \in \mathbb Z_+$ and a function $u \in L^2(\cG_M,m)$ we define $u_M(x) = 2^{Md/2}u(2^{M} x)$. Clearly, $u_M \in L^2(\cG_0,m)$. Also, for $u \in L^2(\cG_0,m)$ let $u_{-M}(x) = 2^{-Md/2} u(2^{-M} x)$, $x \in \cG_M$.
\bigskip

We will need the following scaling properties of Dirichlet forms and principal eigenvalues.

\begin{lem} \label{lem:reflected_dir}
Let $S$ be a complete subordinator with  Laplace exponent $\phi$ given by (\ref{eq:def_phi}). Then the following hold.
\begin{itemize}
\item[(a)] If \textbf{(U1)} is satisfied, then for every $M \in \mathbb Z_+$ and a potential $0 \leq V \in \cK^{X^M}$ we have
$$
\cE_{\phi, V}^M (u,u) = 2^{-Md_w} b \, \cE_{(d_w), \widetilde V}^0 (u_M, u_M), \quad u \in \cF_{\phi,V}^M,
$$
and $
\cF_{(d_w),\widetilde V}^0 = \left\{u \in L^2(\cG_0,m): u_{-M} \in \cF_{\phi,V}^M\right\}$ with $\widetilde V(x) := \frac{2^{Md_w}}{b} \, V(2^M x)$, $x \in \cG_0$. In particular,
$$
\lambda_1^M(\phi, V) = 2^{-Md_w} \, b \, \lambda_1^0(d_w, \widetilde V).
$$

\item[(b)] If \textbf{(U2)} or \textbf{(U3)} is satisfied, then there is a constant $c_{4.4}=c_{4.4}(\phi) \in (0,1]$ such that for every $M \in \mathbb Z_+$ and a potential $0 \leq V \in \cK^{X^M}$, we have
$$
\cE_{\phi, V}^M (u,u) \geq c_{4.4} \, 2^{-M \alpha_1} \cE_{(\alpha_1), \widetilde V}^0 (u_M, u_M), \quad u \in \cF_{\phi,V}^M,
$$
and $
\cF_{(\alpha_1),\widetilde V}^0 \supseteq \left\{u \in L^2(\cG_0,m): u_{-M} \in \cF_{\phi,V}^M\right\}$ with $\widetilde V(x) := \frac{2^{M \alpha_1}}{c_{4.4}} \, V(2^M x)$, $x \in \cG_0$. In particular,
$$
\lambda_1^M(\phi, V) \geq 2^{-M \alpha_1} \, c_{4.4} \, \lambda_1^0(\alpha_1, \widetilde V).
$$
\end{itemize}
\end{lem}

\noindent{\bf Proof.} We only prove (b). The assertion (a) follows directly by definitions of Dirichlet forms and exactly the same arguments.

Assume first that \textbf{(U3)} holds with some $\alpha_1=\alpha_2 \in (0,d_w)$. Let $M \in \mathbb Z_+$ and $0 \leq V \in \cK^{X^M}$. In this case we have $\cF_{\phi, V}^M = \cF_{\phi}^M \cap L^2(\cG_M, V(x)m({\rm d} x))$.
Since $b = 0$, for every $u \in \cF_{\phi,V}^M$, we have, with $J^M$ given by (\ref{eq:j-em}),
$$
\cE_{\phi, V}^M (u,u) = \int_{\cG_M \times \cG_M} (u(x)-u(y))^2 J_{\phi}^M(x,y) m({\rm d}x) m({\rm d}y)
+ \int_{\cG_M} V(x) u^2(x) m({\rm d}x).
$$
By Lemma \ref{lem:comp_levy}, we have $\rho(s) \geq c_{4.1} s^{-1-\alpha_1/d_w}$, $s>0$. Therefore
for every $u \in \cF_{\phi,V}^M$ we get
$$
\int_{\cG_M \times \cG_M} (u(x)-u(y))^2 J_{\phi}^M(x,y) m({\rm d}x) m({\rm d}y) \geq c c_{4.1} \, \cE_{(\alpha_1)}^M (u,u), \quad \mbox{with} \quad c=c(d,\alpha_1),
$$
and, consequently,
\begin{equation}\label{eq:jumps}
\cE_{\phi, V}^M (u,u) \geq c c_{4.1} \, \cE_{(\alpha_1)}^M (u,u) + \int_{\cG_M} V(x) u^2(x) m({\rm d}x).
\end{equation}

We now show that under \textbf{(U2)} the inequality as in (\ref{eq:jumps}) also holds, but an extra step is needed.
Let $u \in \cF_{\phi,V}^M = \cF_{(d_w),V}^M$. Using the estimates from Lemma \ref{lem:comp_levy} and Fubini we will found the lower bound on  $I_M(u):=
\int_{\cG_M \times \cG_M} (u(x)-u(y))^2 J_{\phi}^M(x,y) m({\rm d}x) m({\rm d}y).$ We can write, with any $\delta<1:$
\begin{eqnarray*}
I_M(u)&\geq & c_{4.1} (\delta/2) \int_{\cG_M \times \cG_M} \int_1^\infty (u(x)-u(y))^2 s^{-1-{\alpha_1/d_w}}g^M(s,x,y)\,{\rm d}s\, m({\rm d}x)m({\rm d}y)\\
&=& c_{4.1} (\delta/2) \left( \int_{\cG_M \times \cG_M} \int_0^\infty (u(x)-u(y))^2 s^{-1-{\alpha_1/d_w}}g^M(s,x,y)\,{\rm d}s\, m({\rm d}x)m({\rm d}y)\right.\\
&& \left.- \int_{\cG_M \times \cG_M} \int_0^1 (u(x)-u(y))^2 s^{-1-{\alpha_1/d_w}}g^M(s,x,y)\,{\rm d}s\, m({\rm d}x)m({\rm d}y)\right).
\end{eqnarray*}
In the first of the integrals in the last formula we recognize (up to a constant) the Dirichlet form of the process $X_{(\alpha_1)}^M$, while the other integral is an error term (denoted by $E_M(u)$)  which we will now estimate.  Using Fubini again we write:
\begin{eqnarray}\label{eq:error} E_M(u)&=&
c_{4.1} \delta\int_0^1\left(\frac{1}{2s}\int_{\cG_M\times\cG_M} (u(x)-u(y))^2 g^M(s,x,y)\,m({\rm d}x)m({\rm d}y)\right) s^{-\alpha_1/d_w}\,{\rm d}s.
\end{eqnarray}
For any $s>0$ we have
$$
\frac{1}{2s}\int_{\cG_M\times\cG_M} (u(x)-u(y))^2 g^M(s,x,y)\,m({\rm d}x)m({\rm d}y)\leq  {\mathcal E}^M_{(d_w)}(u,u)
$$
(this is so because the approximating forms increase towards ${\mathcal E}^M_{(d_w)}(u,u)$ as $s\downarrow 0$).
Inserting this bound inside (\ref{eq:error}) and integrating from 0 to 1 we end up with the estimate
\[E_M(u)\leq (c_{4.1} \delta d_w)/ (d_w-\alpha_1)\,\,{\mathcal E}^M_{(d_w)}(u,u).\]
For $\delta= \left(d_w - \alpha_1)/d_w \right) \left( (b/c_{4.1}) \wedge 1\right),$ we have $E_M(u)\leq (b \wedge c_{4.1}) \,\, {\mathcal E}^M_{(d_w)}(u,u),$ therefore we get
\begin{eqnarray}
\cE_{\phi, V}^M (u,u) & = & b \, \cE_{(d_w)}^M (u,u) + \int_{\cG_M \times \cG_M} (u(x)-u(y))^2 J_{\phi}^M(x,y) m({\rm d}x) m({\rm d}y)
+ \int_{\cG_M} V(x) u^2(x) m({\rm d}x) \nonumber \\
& \geq &
b \, \cE_{(d_w)}^M (u,u) - E_M(u) + c c_{4.1}(\delta/2) \cE^M_{(\alpha_1})(u,u)
+ \int_{\cG_M} V(x) u^2(x) m({\rm d}x) \nonumber\\
&\geq & c c_{4.1}(\delta/2)\cE^{M}_{(\alpha_1)}(u,u)+  \int_{\cG_M} V(x) u^2(x) m({\rm d}x) \nonumber.
\end{eqnarray}
This is exactly (\ref{eq:jumps}), with a  smaller constant $c_{4.4}=c c_{4.1}(\delta/2).$ In the sequel, we just write $c_{4.4}$ in either case.

Next, one can directly check using (\ref{eq:scal_gM}) and (\ref{eq:subord}) that $\cE_{(\alpha_1)}^M (u,u) = 2^{-M \alpha_1} \cE_{(\alpha_1)}^0 (u_M,u_M)$.
This way we obtain
$$
\cE_{\phi, V}^M (u,u) \geq 2^{-M \alpha_1} \, c_{4.4} \left( \cE_{(\alpha_1)}^0 (u_M,u_M) + \int_{\cG_0}  \widetilde V (x) u_M^2(x) m({\rm d}x) \right) \nonumber \\
= 2^{-M \alpha_1} \, c_{4.4} \, \cE_{(\alpha_1), \widetilde V}^0 (u_M,u_M),
$$
with $\widetilde V(x):= (2^{M \alpha_1}/c_{4.4}) V(2^M x)$, $x \in \cG_0$. This inequality also implies that
$$
\cF_{(\alpha_1),\widetilde V}^0 \supseteq \left\{u \in L^2(\cG_0,m): u_{-M} \in \cF_{\phi,V}^M\right\}.
$$

To prove the inequality between principal eigenvalues it suffices to use the standard variational formulas for eigenvalues:
$$
\lambda_1^M(\phi, V)  = \inf_{u \in \cF_{\phi, V}^M} \frac{\cE_{\phi, V}^M (u,u)}{\left\|u\right\|^2_{L^2(\cG_M,m)}} \quad \quad \mbox{and} \quad \quad
\lambda_1^0(\alpha_1, \widetilde V) = \inf_{u \in \cF_{(\alpha_1),\widetilde V}^0} \frac{\cE_{(\alpha_1), \widetilde V}^0 (u,u)}{\left\|u\right\|^2_{L^2(\cG_0,m)}}.
$$
Indeed, by the arguments above, for $u \in \cF_{\phi,V}^M$  we have $
u_M \in  \cF_{(\alpha_1),\widetilde V}^0
$ and
$
\cE_{\phi, V}^M (u,u) \geq c_{4.4} \, 2^{-M \alpha_1} \cE_{(\alpha_1), \widetilde V}^0 (u_M, u_M).
$
 Since also $\left\|u\right\|_{L^2(\cG_M,m)} = \left\|u_M\right\|_{L^2(\cG_0,m)},$  for every $u \in \cF_{\phi,V}^M$, we have
$$
\frac{\cE_{\phi, V}^M (u,u)}{\left\|u\right\|^2_{L^2(\cG_M,m)}} \geq c_{4.4} 2^{-M \alpha_1} \frac{\cE_{(\alpha_1), \widetilde V}^0 (u_M,u_M)}{\left\|u_M\right\|^2_{L^2(\cG_0,m)}} \\  \geq  c_{4.4} 2^{-M \alpha_1} \inf_{v \in \cF_{(\alpha_1),\widetilde V}^0} \frac{\cE_{(\alpha_1), \widetilde V}^0 (v,v)}{\left\|v\right\|^2_{L^2(\cG_0,m)}} = c_{4.4} 2^{-M \alpha_1}\, \lambda_1^0(\alpha_1, \widetilde V).
$$
By taking the infimum on the left hand side over all functions $u \in \cF_{\phi,V}^M$, we  get the desired inequality between the principal eigenvalues. The proof is complete.
\hfill$\Box$

\subsubsection{Random Feynman-Kac semigroup and periodization of the Poissonian potential}

Recall that $V$ is called a random Poissonian
potential on $\cG$ if it is given by (\ref{eq:poiss0}). Below we study the process $X^M$ perturbed by the Poissonian potentials $V(x,\omega)$, $x \in \cG_M$, $\omega \in \Omega$, with profiles $W$ satisfying all conditions \textbf{(W1)}--\textbf{(W3)} and restricted to $\left\{(x,y):x,y \in \cG_M \right\}$. As proved in \cite[Proposition 2.1]{bib:KaPP}, under the condition \textbf{(W1)} we have $V(\cdot, \omega) \in \cK_{\loc}^X$ and $V(\cdot, \omega) \in \cK^{X^M}$, $\qpr$-almost surely. The corresponding Feynman-Kac semigroup will be denoted by $(T_t^{\phi, V, M, \omega})_{t \geq 0}$ (resp. $(T_t^{\gamma, V, M,\omega})_{t \geq 0}$ for $X^M_{(\gamma)}$ with $\gamma \in (0,d_w]$). For every $t>0$, the eigenvalues of operators $T_t^{\phi, V, M, \omega}$are given by $\left\{\exp(-t \lambda_n^M(\phi,V,\omega))\right\}_{n \geq 1}$ (resp. $\left\{\exp(-t \lambda_n^M(\gamma,V,\omega))\right\}_{n \geq 1}$), where the random variables $\lambda_n^M(\phi,V,\omega)$ can be ordered as $0 \leq \lambda_1^M(\phi,V,\omega) < \lambda_2^M(\phi,V,\omega) \leq \lambda_3^M(\phi,V, \omega) \leq ... \to \infty$, for $\qpr$-almost all $\omega$.

Our further argument uses some special 'periodization' of the Poissonian potential $V,$ introduced recently in \cite[Def. 3.1]{bib:KaPP}: the family of random fields $(V_M^{*})_{M
\in \mathbb Z_+}$ on $\cG$ given by
\begin{equation} \label{def:Wstar}
V_M^{*}(x,\omega):= \int_{\cG_M} \sum_{y{'} \in \pi_M^{-1}(y)}
W(x,y{'}) \mu^{\omega}({\rm d}y), \quad M \in \mathbb Z_+,
\end{equation}
is called the $M$-periodization of $V$ \textit{in the Sznitman sense}. The same argument as in \cite[Proposition 2.1]{bib:KaPP} yields that  under \textbf{(W1)}, for $\qpr$-almost all $\omega \in \Omega$, one has $V_M^{*}(\cdot,\omega) \in \cK^{X^M}$, for every $M \in \mathbb Z.$

For $t>0$ and $M \in \mathbb Z_+$ we define:
\begin{equation}
L^{N^{*}}_M (t,\omega) = \frac{1}{m(\cG_M)} \int_{\cG_M} p^M(t,x,x)
\ex^{M,t}_{x,x}\left[{\rm e}^{-\int_0^t V_M^{*}(X_s^M,\omega){\rm d}s}\right] m({\rm d}x).
\end{equation}
Our argument in this section essentially relies on the following monotonicity properties.

\begin{lem} \label{lem:monotone}
If  one of the assumptions \textbf{(U1)}-\textbf{(U3)} and all of the assumptions \textbf{(W1)}-\textbf{(W3)} hold, then for any given $t>0 $ we have
\begin{equation}
\mathbb{E}_{\qpr} L^{N^{*}}_{M} (t,\omega) \searrow L(t) \quad \mbox{as} \quad M \to \infty.
\end{equation}
In particular,
\begin{equation} \label{eq:domination_star}
L(t) \leq \mathbb{E}_{\qpr} L^{N^{*}}_{M} (t,\omega), \quad M \in \mathbb Z_{+}, \ t >0.
\end{equation}
\end{lem}

\noindent{\bf Proof.} \cite[Proof of Theorem 3.1]{bib:KaPP} \hfill$\Box$

\medskip

In the sequel we will be mainly working with the following type of rescaled potentials. For a profile $W: \cG \times \cG \to \mathbb R_+$, a random measure $\mu^{\omega}$ with intensity $\nu>0$, and a number $\gamma>0$ we denote
\begin{equation} \label{def:Wstar_new_intensity}
V_{0,M,\gamma}^{*}(x,\omega):= \int_{\cG_0} \sum_{y{'} \in \pi_0^{-1}(y)}
2^{M \gamma} W(2^Mx,2^My') \mu^{M,\omega}({\rm d}y), \quad x \in \cG_0, \quad M \in \mathbb Z_+,
\end{equation}
where $\mu^{M,\omega}$ is the random measure corresponding to the Poisson point process with intensity $2^{Md} \nu$. Clearly, $V_{0,M,\gamma}^{*}$ is the $0$-periodization in the Sznitman sense of the Poissonian potential, which is based on the rescaled profile $2^{M \gamma} W(2^Mx,2^My)$ and the random measure $\mu^{M,\omega}$ with  rescaled intensity.

\subsubsection{Derivation of the upper bound for the short range interaction}

The following upper bound will be the key
point for  our investigations in this subsection.
\begin{lem} \label{lem:upper_basic}
Let $S$ be a complete subordinator with Laplace exponent $\phi$ given by (\ref{eq:def_phi}) and let $V$ be a Poissonian potential with  profile $W$ satisfying the assumptions \textbf{(W1)}-\textbf{(W3)}. The following hold.
\begin{itemize}
\item[(a)] Under the assumption \textbf{(U1)}, there exists a constant $c_{4.5}>0$ such that for every $t>1$ and every number $M \in \mathbb Z_+$ such that $M \leq \frac{\log_2 (t/\nu)}{d+d_w}$ we have
\begin{equation}
\mathbb{E}_{\qpr} [L^{N^{*}}_M (t,\omega)] \leq c_{4.5} \, \mathbb{E}_{\qpr} \exp \left[- b\left(1-\frac{1}{t}\right) \, \nu^{\frac{d_w}{d+d_w}} t^{\frac{d}{d+d_w}} \, \lambda_1^0 \left(d_w, V_{0,M,d_w}^{*}, \omega\right) \right],
\end{equation}
where the potential $V_{0,M,d_w}^{*}$ is given by (\ref{def:Wstar_new_intensity}).

\item[(b)] Under the assumption \textbf{(U2)} or \textbf{(U3)}, there exists a constant $c_{4.6}>0$ such that for every $t>1$ and every number $M \in \mathbb Z_+$ such that $M \leq \frac{\log_2 (t/\nu)}{d+\alpha_1}$ we have
\begin{equation}
\mathbb{E}_{\qpr} L^{N^{*}}_M (t,\omega)  \leq c_{4.6} \, \mathbb{E}_{\qpr} \exp \left[-  c_{4.4}\left(1-\frac{1}{t}\right) \, \nu^{\frac{\alpha_1}{d+\alpha_1}} t^{\frac{d}{d+\alpha_1}} \, \lambda_1^0 \left(\alpha_1, V_{0,M,\alpha_1}^{*} , \omega\right) \right],
\end{equation}
where the potential $V_{0,M,\alpha_1}^{*}$ is given by (\ref{def:Wstar_new_intensity}).
\end{itemize}
\end{lem}

\noindent{\bf Proof.} We only prove (b). The proof of (a) requires exactly the same argument and is even easier. Let $\phi$ satisfy \textbf{(U2)} or \textbf{(U3)} and let $V$ be a Poissonian potential with  profile $W$ as in the assumptions. Fix arbitrary $t>1$ and $M \in \mathbb Z_+$ such that $M \leq \frac{\log_2 (t/\nu)}{d+\alpha_1}$. By Fubini, for all such $t$ and $M$, we get
\begin{eqnarray} \label{eq:eq1}
\mathbb{E}_{\qpr} L^{N^{*}}_M (t,\omega) = \frac{1}{m(\cG_M)} \int_{\cG_M} p^M(t,x,x)
\ex^{M,t}_{x,x}\left[\mathbb{E}_{\qpr}\left[{\rm e}^{-\int_0^t V_M^{*}(X^M_s,\omega){\rm d}s}\right]\right] m({\rm d}x).
\end{eqnarray}
Observe now that by the exponential formula (\ref{eq:exp-formula}) and the scaling properties of the measure $m$, for every measurable and nonnegative function $f$ on $\cG$ we have
\begin{eqnarray*}
\mathbb{E}_{\qpr} e^{-\int_{\cG} f(y) \mu^{\omega}({\rm d}y)} & = & \exp \left(-\nu \int_{\cG} \left(1-e^{-f(y)}\right) m ({\rm d} y) \right) \\
& = & \exp \left(-2^{Md} \nu \int_{\cG} \left(1-e^{-f(2^M y)}\right) m ({\rm d} y) \right) = \mathbb{E}_{\qpr} e^{-\int_{\cG} f^M(y) \mu^{M, \omega}({\rm d}y)},
\end{eqnarray*}
where $f^M(y) = f(2^M y)$, $y \in \cG$, and $\mu^{M, \omega}$ is the random measure corresponding to the Poisson point process with rescaled intensity $2^{Md} \nu$. Applying this observation to the functions
$$\cG\ni y\mapsto
f_{w}(y):=  \1_{\cG_M}(y) \cdot \sum_{y{'} \in \pi_M^{-1}(\pi_M(y))}
\int_0^t W(X^M_s(w) ,y{'}) ds,
$$ we obtain that
for every $x \in \cG_M$ and $\pr^{M,t}_{x,x}$-almost all $w$
\begin{eqnarray} \label{eq:eq2}
\mathbb{E}_{\qpr} \left[{\rm e}^{-\int_0^t V_M^{*}(X_s^M,t){\rm d}s}\right] = \mathbb{E}_{\qpr} \left[ {\rm e}^{-\int_0^t V_{0}^{M*}(X_s^M,\omega){\rm d}s}\right],
\end{eqnarray}
where
$$
V_{0}^{M*}(x,\omega):= \int_{\cG_0}  \sum_{y{'} \in \pi_M^{-1}(2^M y)}
W(x ,y{'}) \mu^{M, \omega}({\rm d} y) = \int_{\cG_0}  \sum_{y{'} \in \pi_0^{-1}(y)}
W(x , 2^M y{'}) \mu^{M, \omega}({\rm d} y).
$$
Inserting (\ref{eq:eq2}) to (\ref{eq:eq1}) and using the bridge kernel representation \cite[(2.28)]{bib:KaPP}, we thus get
\begin{eqnarray*}
\mathbb{E}_{\qpr} L^{N^{*}}_M (t,\omega) &=  & \frac{1}{m(\cG_M)} \int_{\cG_M} p^M(t,x,x)
\ex^{M,t}_{x,x}\left[\mathbb{E}_{\qpr}\left[{\rm e}^{-\int_0^t V_{0}^{M*}(X_s^M,\omega){\rm d} s}\right]\right] m({\rm d}x) \\
& = &\frac{1}{m(\cG_M)}
\mathbb{E}_{\qpr}  \mbox{Tr} \, T_t^{\phi,V_{0}^{M*}, M, \omega}= \frac{1}{m(\cG_M)} \mathbb{E}_{\qpr}  \sum_{n=1}^{\infty} e^{-t \, \lambda_n^M \left(\phi, V_{0}^{M*} , \omega \right)} \nonumber \\
& = &  \frac{1}{m(\cG_M)}\mathbb{E}_{\qpr} \sum_{n=1}^{\infty} e^{-(t-1) \, \lambda_n^M \left(\phi, V_{0}^{M*} , \omega \right)} \, e^{-\lambda_n^M \left(\phi, V_{0}^{M*} , \omega \right)}.
\end{eqnarray*}
Since for $\qpr$-almost all $\omega \in \Omega$ we have $0 \leq \lambda_1^M \left(\phi, V_{0}^{M*} , \omega \right) < \lambda_2^M \left(\phi, V_{0}^{M*} , \omega \right) \leq \lambda_3^M \left(\phi, V_{0}^{M*} , \omega \right) \leq ...$, it follows that
\begin{eqnarray*}
\mathbb{E}_{\qpr} L^{N^{*}}_M (t,\omega) & \leq & \mathbb{E}_{\qpr} e^{-(t-1) \lambda_1^M \left(\phi, V_{0}^{M*} , \omega \right)} \cdot \frac{1}{m(\cG_M)} \mbox{Tr} \, T_1^{\phi,V_{0}^{M*}, M, \omega} \\ & \leq & \mathbb{E}_{\qpr} e^{-(t-1) \lambda_1^M \left(\phi, V_{0}^{M*} , \omega \right)} \cdot \frac{1}{m(\cG_M)} \int_{\cG_M} p^M(1,x,x) m({\rm d}x).
\end{eqnarray*}
Moreover, by Lemma \ref{lem:pMbound} we have $p^M(1,x,x) \leq 3 c_{4.3}$, for every $x \in \cG_M$ and $M \in \mathbb Z_+$. Thus, we get
$$
\mathbb{E}_{\qpr} L^{N^{*}}_M (t,\omega) \leq 3 c_{4.3} \mathbb{E}_{\qpr} e^{-(t-1) \lambda_1^M \left(\phi, V_{0}^{M*} , \omega \right)}.
$$
 Since $V_{0,M,\alpha_1}^{*}(x) = 2^{M \alpha_1} V_{0}^{M*}(2^M x)$, $x \in \cG_0$ (recall that $V_{0,M,\alpha_1}^{*}$ is given by (\ref{def:Wstar_new_intensity})), we derive from Lemma \ref{lem:reflected_dir} (b) the inequality
$$
\lambda_1^M \left(\phi, V_{0}^{M*} , \omega \right) \geq c_{4.4} 2^{-M\alpha_1} \lambda_1^0 \left(\alpha_1, V_{0,M,\alpha_1}^{*}, \omega \right), \quad M \in \mathbb Z_+,
$$
which holds for $\qpr$-almost all $\omega$. In consequence,
$$
\mathbb{E}_{\qpr} L^{N^{*}}_M (t,\omega) \leq 3 c_{4.3} \mathbb{E}_{\qpr} e^{-c_{4.4} (t-1) 2^{-M\alpha_1} \lambda_1^0 \left(\alpha_1, V_{0,M,\alpha_1}^{*} , \omega \right)}.
$$
Therefore, for every $t>1$ and $M \in \mathbb Z_+$ such that $2^M \leq (t/\nu)^{\frac{1}{d+\alpha_1}}$ we finally get
$$
\mathbb{E}_{\qpr} L^{N^{*}}_M (t,\omega) \leq 3 c_{4.3} \mathbb{E}_{\qpr} \exp \left[-  c_{4.4}\left(1-\frac{1}{t}\right) \, \nu^{\frac{\alpha_1}{d+\alpha_1}} t^{\frac{d}{d+\alpha_1}} \, \lambda_1^0 \left(\alpha_1, V_{0,M,\alpha_1}^{*}, \omega\right) \right].
$$ The proof is complete.

\hfill$\Box$

\medskip

Under the following additional assumption on the profile $W$:
\begin{description}
\item[(W4)] {There exist constants $a_0, A>0$ such that
\begin{eqnarray*} W(x,y)\geq A \mbox{ when } d(x,y)\leq a_0
\end{eqnarray*}}
\end{description}
we prove the following theorem.

\begin{theo}\label{gorne-sbm}
Let $X$ be a subordinate Brownian motion in $\cG$ via the subordinator $S$ with Laplace exponent $\phi$ of the form (\ref{eq:def_phi}) and let $V$ be a Poissonian potential with the profile $W$ such that the assumptions \textbf{(W1)}-\textbf{(W4)} are satisfied. Then there exists $D_1>0$ such that the  following hold.

\begin{itemize}
\item[(a)]
Under the assumption \textbf{(U1)}:
\begin{equation}\label{gor1}
\limsup_{t\to \infty}\frac{\log
L(t)}{t^{\frac{d}{d+d_w}}}\leq
-D_1 \, \nu^{\frac{d_w}{d+d_w}}.
\end{equation}
\item[(b)] Under the assumption \textbf{(U2)} or \textbf{(U3)}:
\begin{equation}\label{gor2}
\limsup_{t\to \infty}\frac{\log
L(t)}{t^{\frac{d}{d+\alpha_1}}}\leq
-D_1 \, \nu^{\frac{\alpha_1}{d+\alpha_1}}.
\end{equation}

\end{itemize}
\end{theo}

\noindent{\bf Proof.}  In both cases (a), (b) we use  Sznitman's theorem from the Appendix, in either its diffusion version \cite[Theorem 1.4]{bib:Szn1}
or the non-diffusion version \cite[Theorem 1]{bib:Kow}, adapted to the potential case. As the statements of both these theorems are nearly identical (they pertain to either $d_w$ or $\alpha_1 \in (0,d_w)$), we will write $`\gamma$' for  $d_w$ or $\alpha_1\in (0,d_w),$ depending on the context.

Let now
\begin{eqnarray} \label{eq:choiceM}
M=M(t)=\left[\frac{\log_2(t/\nu)}{(d+\gamma)}\right],&\mbox{ i.e. } & 2^{M}\leq \left(\frac{t}{\nu} \right)^{1/(d+\gamma)}< 2^{M+1},
\end{eqnarray}
and write $\epsilon = 2^{-M}.$
By Lemmas \ref{lem:monotone} and \ref{lem:upper_basic},  for every $t > 1$ and  $M=M(t)$ given by (\ref{eq:choiceM}) we have
$$
L(t) \leq \mathbb{E}_{\qpr} L^{N^{*}}_M (t,\omega) \leq c_{4.6} \, \mathbb{E_Q}\exp\left[-c_{4.4}\left(1-\frac{1}{t}\right)
\nu^{\frac{\gamma}{d+\gamma}}
t^{\frac{d}{d+\gamma}}\lambda_1^0(\gamma,V^*_{0,M,\gamma},\omega) \right],
$$
with $b$ and $c_{4.5}$ replacing $c_{4.4}$ and $c_{4.6}$ in  case (a). Hence, it is enough to estimate the expectation on the right hand side of the formula above, independently of $M.$

{Fix a number  $a > a_0$ and denote $W_a(x,y) = W(x,y) \cdot \textbf{1}_{d(x,y) \leq a}$, $x,y \in \cG$. } The number $a$ will not vary throughout the proof. The periodized potential $V^*_{0,M,\gamma}(x),$ $x\in \mathcal G_0,$ satisfies
\[V^*_{0,M,\gamma}(x,\omega)\geq \int_{\mathcal G_0} 2^{M\gamma}{W_a(2^Mx,2^My)}\,\mu^{M,\omega}({\rm d}y)=: \widetilde V^*_{0,M,\gamma}(x,\omega),\]
where now $ \mu^{M,\omega}$ comes from the rescaled cloud on $\mathcal G_0$ with intensity $\tilde\nu=2^{Md}\nu,$ whose law will still be denoted by ${\mathbb Q}.$
The new profile $\mathcal G_0\times \mathcal G_0\ni (x,y)\mapsto {W_{a,M}(x,y)= 2^{M\gamma}W_a(2^Mx, 2^My)}$ has range $a2^{-M}=a\epsilon,$ and its values are bigger than $2^{M\gamma}A$ when $d(x,y)\leq a_0\epsilon.$

Further assume that $b=2^{\kappa}>a,$ with $\kappa\in\mathbb Z,$ and let
$K,\delta>0$ be given.

$\mathbb Q-$almost surely, there is a finite number of Poisson points in $\mathcal G_0.$ From now on we will be working with a  fixed
configuration $\omega=(x_1,...,x_N)\subset \mathcal G_0$ of Poisson points. We divide them into `good' and 'bad' points
according to Definition \ref{defi:good-bad}, and remove the closed balls $B(x_i,b\epsilon)$ with centers at good points from the state-space. We are left with the set
\[\Theta_{b,M}= \mathcal G_0 \setminus \bigcup_{x_i-{\rm good}}\overline{B}(x_i, b\epsilon),\]
and we let the process $X_{(\gamma)}^0$ evolve in this set, being killed  when it enters
one of the balls $\overline{B}(x_i, b\epsilon),$ $x_i-$good.
Let $\lambda_1^0(\gamma, \Theta_{b,M}, \omega)$ be the principal eigenvalue of the generator of this process.

The assumptions of Theorem \ref{compare} are fulfilled
(we postpone their verification until after the proof; see Subsection \ref{subsec:ver}),
and so
 there exists $\epsilon_0>0,$ depending on the process, the potential $W,$ and the numbers $a,b,K,\delta,\gamma$ (not on $M$) such that when $\epsilon<\epsilon_0,$ then
\begin{eqnarray}\label{eq:lambdas}
\lambda_1^0(\gamma, \Theta_{b,M}, \omega) \wedge K&\leq& \lambda_1^0(\gamma, \widetilde V^*_{0,M,\gamma},\omega)\wedge K+\delta \nonumber \\
&\leq &\lambda_1^0(\gamma,  V^*_{0,M,\gamma},\omega)\wedge K+\delta
\end{eqnarray}
(the last inequality follows from the inequality  $V^*_{0,M,\gamma}\geq \widetilde V^*_{0,M,\gamma}$, combined with the
variational definition of the principal eigenvalue).
In particular, since $\epsilon=2^{-M},$ there exists $M_0$ such that for $M>M_0$  the relation (\ref{eq:lambdas}) holds. The way $M=M(t)$ was defined (see (\ref{eq:choiceM})) gives that there exists {$t_0 \geq 1$} such that it holds for $t>t_0.$

The conclusion of the proof is much alike the conclusion of \cite[Theorem 1.7]{bib:Szn1} or \cite[Lemma 9]{bib:KPP-PTRF}.
Let $M>M_0$ (equivalently: $t>t_0$).
We chop the sides of  the triangle $\mathcal G_0$ into $(b\epsilon)^{-1}=2^{M-\kappa}$ parts, which yields $N(b,M)=2^{(M-\kappa)d}=(b\epsilon)^{-d}$ small gasket triangles of
sidelength $b\epsilon.$
Now, instead of removing balls $\overline B(x_i,b\epsilon)$ from the state-space, we remove those closed small triangles that received some
 (good) Poisson points. More precisely,
let
$A_{b,M}$ be the union of those small triangles that received some
Poisson points, and $\widehat{A}_{b,M}$ -- of those triangles that
received some good Poisson points. We set
\begin{eqnarray*}
U_{b,M}=\mathcal G_0\setminus A_{b,M},&& \widehat{U}_{b,M}=\mathcal G_0\setminus \widehat {A}_{b,M}.
\end{eqnarray*}

As the diameter
of each of the  triangles
removed equals to $b\epsilon,$ we have $\Theta_{b,M}\subset \widehat{U}_{b,M},$ and consequently $\lambda_1^0(\gamma, \Theta_{b,M}, \omega)\geq \lambda_1^0(\gamma,\widehat {U}_{b,M},\omega),$ where $\lambda_1^0(\gamma, \widehat {U}_{b,M},\omega)$ is the principal eigenvalue of the process that is killed upon exiting $\widehat {U}_{b,M}.$

Altogether, for any given configuration $\omega,$
given $K,\delta, b,$ we have
\begin{eqnarray*}
\lambda_1^0(\gamma,V^*_{0,M,\gamma},\omega) &\geq & \lambda_1^0(\gamma,V^*_{0,M,\gamma},\omega)\wedge K
\geq  \lambda_1^0(\gamma, \Theta_{b,M}, \omega)\wedge K-\delta\\
&\geq& \lambda_1^0(\gamma,  \widehat {U}_{b,M},\omega)\wedge K-\delta.
\end{eqnarray*}

Denoting by $\mathcal U_{b,M}$ the collection of all possible configurations of the sets $U_{b,M}$ and  $\widehat {U}_{b,M}$
and noting that $\# \mathcal U_{b,M}= 2^{N(b,M)},$ we can proceed as follows (taking the precisely chosen $M=M(t)$ and $t>t_0$):
\begin{eqnarray*}&& \mathbb{E_Q}\exp\left[-c_{4.4}\left(1-\frac{1}{t}\right)
\nu^{\frac{\gamma}{d+\gamma}}
t^{\frac{d}{d+\gamma}}\lambda_1^0(\gamma,V^*_{0,M,\gamma},\omega) \right]\\ &\leq& \mathbb{E_Q}\exp\left[-c_{4.4}\left(1-\frac{1}{t}\right)
\nu^{\frac{\gamma}{d+\gamma}}
t^{\frac{d}{d+\gamma}}(\lambda_1^0(\gamma, \widehat {U}_{b,M},\omega)\wedge K-\delta) \right]\\
&\leq & \sum_{U,\widehat{U}\in\mathcal U_{b,M}}\mathbb{E_Q}\left[\exp\left(-c_{4.4}\left(1-\frac{1}{t}\right)
\nu^{\frac{\gamma}{d+\gamma}}
t^{\frac{d}{d+\gamma}}(\lambda_1^0(\gamma, \widehat{U}_{b,M},\omega)\wedge K-\delta) \right) \mathbf 1\{U_{b,M}= U, \widehat {U}_{b,M}=\widehat{U}\}\right].
\end{eqnarray*}

For any $A\in\mathcal B(\mathcal G_0)$ we have $\mathbb Q[\mathcal N(A)=0]={\rm e}^{-\tilde \nu m(A)}.$
One knows \cite[Lemma 1.3]{bib:Szn1} that \linebreak $m\left(\bigcup_{x_i-{\rm bad}}\overline B(x_i, b\epsilon)\right)\leq \delta,$ therefore also
$m(A_{b,M}\setminus \widehat{A}_{b,M})\leq \delta$ and
$m(\widehat{U}_{b,M})\leq m( U_{b,M})+\delta.$ Therefore the estimate continues as
$$
= \sum_{(U,\widehat{U})\in\mathcal A_{b,M}}\exp\left(-c_{4.4}\left(1-\frac{1}{t}\right)
\nu^{\frac{\gamma}{d+\gamma}}
t^{\frac{d}{d+\gamma}}(\lambda_1^0(\gamma, \widehat {U})\wedge K-\delta) \right)
\mathbb Q[U_{b,M}= U, \widehat{U}_{b,M}=\widehat{U}] \ =: I,
$$
where $\mathcal A_{b,M}\subset \mathcal U_{b,M}\times \mathcal U_{b,M}$ consists of those pairs $(U,\widehat{U})$ for which $U\subset\widehat{U},$ $m(\widehat U)\leq m(U)+\delta,$ and $\lambda_1^0(\gamma, \widehat {U})$ is the principal eigenvalue of the generator of the process $X_{(\gamma)}^0$ killed upon exiting the open set $\widehat{U}$.
In particular, for $(U,\widehat{U})\in\mathcal A_{b,M}$ one has
\[\mathbb Q[U_{b,M}=U, \widehat{U}_{b,M}=\widehat{U}] \leq {\rm e}^{-\tilde \nu (m(\widehat U)-\delta)}.\]
Since
\[\#\mathcal U_{b,M}= 2^{N(b,M)}\leq 2^{b^{-d}\left(\frac{t}{\nu}\right)^{d/(d+\gamma)}}\]
and $2^{-d}\nu^{\frac{\gamma}{d+\gamma}}t^{\frac{d}{d+\gamma}} \leq \tilde\nu=2^{Md}\nu \leq \nu^{\frac{\gamma}{d+\gamma}}t^{\frac{d}{d+\gamma}},$ we get
\begin{eqnarray*}I
&\leq&   \sum_{U,\widehat U\in\mathcal U_{b,M}}\exp\left(-\left[c_{4.4}\left(1-\frac{1}{t}\right)
\nu^{\frac{\gamma}{d+\gamma}}
t^{\frac{d}{d+\gamma}}(\lambda_1^0(\gamma,  \widehat{U})\wedge K-\delta)\right]-\tilde \nu(m(\widehat{U}) -\delta) \right)\\
&\leq & (2^{\# \mathcal U_{b,M}})^2 \exp\left(-c_{4.4}\left(1-\frac{1}{t}\right) \nu^{\frac{\gamma}{d+\gamma}}t^{\frac{d}{d+\gamma}}\inf_{U\in\mathcal U_{b,M}} [(\lambda_1^0(\gamma,U)\wedge K-\delta)+\frac{2^{-d}t}{c_{4.4}(t-1)} (m(U)- 2^{d}\delta)]\right)\\
&\leq & \exp\left[\frac{2\log 2}{b^d} \left(\frac{t}{\nu}\right)^{\frac{d}{d+\gamma}}  \!\!\!\! - c_{4.4}\left(1-\frac{1}{t}\right) \nu^{\frac{\gamma}{d+\gamma}}t^{\frac{d}{d+\gamma}}\inf_{U\in\mathcal U_0} [(\lambda_1^0(\gamma,U)\wedge K-\delta)+\frac{2^{-d}t}{c_{4.4}(t-1)}  (m(U)- 2^{d}\delta)]\right],
\end{eqnarray*}
where by $\mathcal U_0$ we have denoted the collection of all open subsets of $\mathcal G_0.$ This bound is valid for $t>t_0.$ In particular,
\begin{eqnarray*}\limsup_{t\to\infty}\frac{\log L(t)}{t^{\frac{d}{d+\gamma}}}&\leq&
\frac{2\ln 2}{b^d\nu^{\frac{d}{d+\gamma}}}-c_{4.4}\nu^{\frac{\gamma}{d+\gamma}}
\inf_{U\in\mathcal U_0}[(\lambda_1^0(\gamma,U)\wedge K-\delta)+ (2^dc_{4.4})^{-1}  (m(U)- 2^{d}\delta)]\\
&\leq & \frac{2\ln 2}{b^d\nu^{\frac{d}{d+\gamma}}}-c_{4.4}\nu^{\frac{\gamma}{d+\gamma}}
\left(\inf_{U\in\mathcal U_0}[\lambda_1^0(\gamma,U)+ (2^dc_{4.4})^{-1} m(U)]\wedge K-\delta(1+ c_{4.4}^{-1})\right).
\end{eqnarray*}
The left-hand side does not depend on $b,K,\delta$ --  by passing to the limit $b\to\infty,$ $\delta\to 0,$ $K\to\infty $ on the right-hand side
we get
\[\limsup_{t\to\infty}\frac{\log L(t)}{t^{\frac{d}{d+\gamma}}}\leq
-c_{4.4} \nu^{\frac{\gamma}{d+\gamma}}\inf_{U\in\mathcal U_0}[\lambda_1^0(\gamma,U)+ (2^dc_{4.4})^{-1} m(U)]=: -D_1\nu^{\frac{\gamma}{d+\gamma}}.\]
To conclude, we need to check that $D_1>0.$
We can write, for any open $U\subset \mathcal G_0,$
\begin{eqnarray*}{\rm e}^{-{t\lambda_1^0(\gamma,U)}} &\leq & \mbox{Tr}\, P_t^{\gamma,0,U} \leq \int_{U} p^{\gamma, 0}(t,x,x)\,{\rm d}m(x)
\end{eqnarray*}
where $p^{\gamma,0}$ is the transition density of the process $X^0_{(\gamma)}$ on $\mathcal G_0,$ and  $P_t^{\gamma,0,U}$ is the semigroup of this process killed upon exiting $U.$ In particular, for $t>1, $ from (\ref{eq:pMbound1}) with $\alpha_1=\alpha_2=\beta=\gamma$ we get $\sup_{x\in\mathcal G_0}p^{\gamma,0}(t,x,x)\leq 3 c_{4.3},$ so that
\[\lambda_1^0(\gamma,U) + (2^dc_{4.4})^{-1} m(U) \geq - \frac{\log (3 c_{4.3}m(U))}{t}+ (2^dc_{4.4})^{-1} m(U).\]
As for $t>2^dc_{4.4}$ one has\[\inf_{x\in(0,1)}\left[ -\frac{\log (3 c_{4.3}x)}{t} + (2^dc_{4.4})^{-1} x\right]= \frac{1}{t}\left(1-\log\textstyle\frac{ 2^dc_{4.4} \cdot 3 c_{4.3}}{ t}\right),\]
 picking any $t> 2^dc_{4.4} \max(1, 3 c_{4.3}e^{-1})$ we achieve
the desired statement that $D_1>0.$ \hfill$\Box$

\medskip

We also prove the matching bounds for the Feynman-Kac
functionals.

\begin{theo}\label{th:funct-upper}
Let $X$ be a subordinate Brownian motion in $\cG$ via the subordinator $S$ with Laplace exponent $\phi$ of the form (\ref{eq:def_phi}) and let $V$ be a Poissonian potential with the profile $W$ such that the assumptions \textbf{(W1)}-\textbf{(W4)} are satisfied. Then the following hold (with $D_1$ same as above).

\begin{itemize}
\item[(a)]
Under the assumption \textbf{(U1)}
for any $x\in\mathcal G$
\begin{equation}\label{gor1-FKF}
\limsup_{t\to\infty} \frac{\log \mathbb{E_Q}\mathbf E_x[{\rm e}^{-\int_0^t
V(X_s,\omega)\,{\rm d}s}]}{t^{\frac{d}{d+d_w}}}\leq
-D_1\, \nu^{\frac{d_w}{d+d_w}}.
\end{equation}
\item[(b)] Under the assumption \textbf{(U2)}
or \textbf{(U3)}
 for any $x\in \mathcal G$
\begin{equation}\label{gor2-FKF}
\limsup_{t\to\infty} \frac{\log \mathbb{E_Q}\ex_x[{\rm e}^{-\int_0^t
V(X_s,\omega)\,{\rm d}s}]}{t^{\frac{d}{d+\alpha}}}\leq
-D_1\, \nu^{\frac{{\alpha_1}}{d+{\alpha_1}}}.
\end{equation}
\end{itemize}
\end{theo}

\noindent{\bf Proof.}
By the same argument as in \cite[Ineq. (3.8)]{bib:KaPP}, the $M$-periodicity of the potential $V^{*}_M$ and the definition of the measure $\pr^M$, for $M$ so large that $x\in\mathcal G_M$ and for $t>1$, we get
\begin{eqnarray*}\mathbb{E_Q}\mathbf E_x\left[{\rm e}^{-\int_0^t V(X_s,\omega)\,{\rm d}s}\right]
& \leq & \mathbb{E_Q}\mathbf E_x\left[{\rm e}^{-\int_0^t V^{*}_M(X_s,\omega)\,{\rm d}s}\right] = \mathbb{E_Q}\mathbf E^M_x\left[{\rm e}^{-\int_0^t V^{*}_M(X^M_s,\omega)\,{\rm d}s}\right] \\ & = & \int_{\mathcal G_M \times \mathcal G_M} u^{M,\omega}(1,x,z) u^{M,\omega}(t-1,z,y) {\rm d}m(z){\rm d}m(y),
\end{eqnarray*}
where
\[u^{M,\omega}(t,x,y)= p^M(t,x,y)\mathbf E_{x,y}^{M,t}\left[{\rm e}^{-\int_0^t V^*_M(X^M_s,\omega){\rm d}s}\right]\]
is the kernel of the operator $T_t^{\phi,V^{*}_M,M,\omega}$ (see \cite[(2.28)]{bib:KaPP}). From Lemma \ref{lem:pMbound} we have $u^{M,\omega}(1,x,z) \leq p^M(1,x,z) \leq 3c_{4.3}$, for every $x,z \in \cG_M$. Thus, by the same estimate for the kernel $u^{M,\omega}(t-1,z,y)$ as in \cite[the last two lines on p. 235]{bib:Szn1}, we may conclude that
$$
\mathbb{E_Q}\mathbf E_x\left[{\rm e}^{-\int_0^t V(X_s,\omega)\,{\rm d}s}\right] \leq3c_{4.3}m(\mathcal G_M) \mathbb {E_Q}
\int_{\mathcal G_M} p^M(t-1,x,x)\mathbf E_{x,x}^{M,t-1}\left[{\rm e}^{-\int_0^{t-1} V^*_M(X^M_s,\omega){\rm d}s}\right]{\rm d}m(x).
$$
In the right-hand side we recognize the expression $3c_{4.3} [m(\mathcal G_M)]^2\mathbb {E_Q}[L^{N^*}_M(t-1,\omega)],$ which has been already
estimated in Lemma \ref{lem:upper_basic}. Therefore, starting from (\ref{eq:choiceM}) and proceeding exactly in the same way as in the proof of Theorem \ref{gorne-sbm}, we finally get the desired inequality with $\gamma=d_w$ for (a) and $\gamma = \alpha$ for (b) respectively, with the same constant $D_1$ as before.
\hfill$\Box$

\subsection{Verification of the assumptions of Sznitman's theorem} \label{subsec:ver}
We apply the Sznitman's theorem (Theorem \ref{compare} below) in the following setting:
\begin{description}
\item[$*$]$\mathcal T= \mathcal G_0,$ the metric $d$ is the shortest path distance on $\mathcal G_0,$ $m$ is the Hausdorff measure on $\mathcal G_0$ in dimension $d=\frac{\log 3}{\log 2}$ normalized to have $m(\mathcal G_0)=1,$ which is a doubling measure;
 \item[$*$] the Markov process in question is $X^0_{(\gamma)}$ on $\mathcal G_0$ for $\gamma \in (0,d_w]$ (the reflected jump stable process or the reflected Brownian motion);
\item[$*$] for $x,y\in\mathcal G_0$ the potential profile is given by $W_M(x,y):= 2^{M\gamma} W(2^Mx,2^My),$  where the profile $W$ is of finite range $a>0$ (i.e. $W(x,y)=0$ when $d(x,y)>a$) and satisfies {\bf (W1)} -- {\bf (W4)}. Assume that $a\geq a_0,$ where $a_0$ is the constant from {\bf (W4)}.
    The range of $W_M$ is equal to $a2^{-M};$ we denote $2^{-M}=\epsilon,$ and we will also write $W_\epsilon$ for $W_M.$
\end{description}
All the required regularity assumptions (see Subsection \ref{subsec:ass}) except for {\bf (P3)} were established in \cite{bib:KPP-PTRF, bib:Szn1} for the reflected Brownian motion and in \cite{bib:Kow-KPP} for reflected jump stable processes. These papers were concerned with processes
on $\mathcal G$ evolving among killing Poissonian obstacles. We now verify the remaining condition {\bf (P3)}, which is needed in our case.

\begin{prop}\label{prop:P3} Let $\gamma\in(0,d_w]$ and
let $X^0=X^0_{(\gamma)}$ be the reflected $\gamma-$stable process on $\mathcal G_0$ (not excluding the case $\gamma=d_w$). Assume that the potential profile $W$ satisfies the condition {\bf (W4)}. Then there exists
 constants $c_{4.7}=c_{4.7}(a_0,A,b,\gamma)>0$ and $\tau_0=\tau_0(a_0,b,\gamma)>0$ such that for any $x\in \mathcal G_0,$ $\epsilon=2^{-M}>0,$ and $y\in\mathcal G_0$ with $d(x,y)\leq b\epsilon$  one has
 \[\ex_x^0\left[{\rm e}^{-\int_0^{\tau_0\epsilon^\gamma/2} W_\epsilon(X_s^0,y)\,{\rm d}s}\right]\leq 1-2c_{4.7}, \ \ \ \ \ \quad \mbox{with} \quad W_\epsilon(x,y)= \epsilon^{-\gamma} W(x/\epsilon, y/\epsilon) .\]
 \end{prop}

 Note that the constant $c_{4.7}$ does not depend on $\epsilon,$ which is
decisive for the proof of the upper bound theorem ($c_{4.7}$ plays the role of the constant $c_1$ in assumptions \textbf{(P2)}-\textbf{(P3)} in Subsection \ref{subsec:ass}).

\

\noindent {\bf Proof.}  Suppose $M\geq 0$ and $x,y\in \mathcal G_0$ are as in the assumptions.
Reflected processes on the gaskets $\mathcal G_M$ allow for discrete scaling. Namely, for $x\in\cG_0,$
the processes on $\cG_M:$ $(2^{M}X^0_s)_{s\geq 0}$  under $\mathbf P_x^0$ and   $(X^M_{2^{M\gamma}s})_{s\geq 0}$ under $\mathbf P_{2^Mx}^M$ are equal in law. Denote $\tilde x= 2^Mx (=\epsilon^{-1}x)$ and $\tilde y= 2^M y$ so that for any $t>0$
\begin{eqnarray*}\mathbf E_x^0\left[{\rm e}^{-\int_0^{t\epsilon^\gamma}W_\epsilon (X_s^0,y)\,{\rm d}s}\right] &=& \mathbf E_{2^{-M}\tilde x}^0 \left[{\rm e}^{-\int_0^{t 2^{-M\gamma}} 2^{M\gamma} W(2^M X^0_{s}, 2^My)\,{\rm d}s}\right]\\
&=& \mathbf E_{\tilde x}^M\left[{\rm e}^{-\int_0^{t 2^{-M\gamma}} 2^{M\gamma} W(X^M_{2^{M\gamma}s}, \tilde y)\,{\rm d}s}\right] = \mathbf E_{\tilde x}^M\left[{\rm e}^{-\int_0^{t} W(X^M_s, \tilde y)\,{\rm d}s}\right],
\end{eqnarray*}
where now  $\tilde x, \tilde y \in \mathcal G_M$ and $d(\tilde x, \tilde y)\leq b.$

For $r>0$ and fixed $\tilde y\in \mathcal G_M,$ denote $T_r=\inf\{t\geq 0: X_t^M\in B(\tilde y, r)\}.$
Then  we can write, in the third line using the strong Markov property at the stopping time $T_{a_0/2}$:
\begin{eqnarray}
\mathbf E_{\tilde x}^M\left[{\rm e}^{-\int_0^t W(X_s^M,\tilde y)\,{\rm d}s}\right] &= & \mathbf E_{\tilde x}^M\left[{\rm e}^{-\int_0^t W(X_s^M,\tilde y)\,{\rm d}s}\mathbf 1\{T_{a_0/2}<\textstyle\frac{t}{2}\}\right]
+ \mathbf E_{\tilde x}^M\left[{\rm e}^{-\int_0^t W(X_s^M,\tilde y)\,{\rm d}s}\mathbf 1\{T_{a_0/2}\geq \textstyle\frac{t}{2}\}\right]\nonumber\\[1mm]
&\leq & \mathbf E_{\tilde x}^M\left[{\rm e}^{-\int_0^t W(X_s^M,\tilde y)\,{\rm d}s}\mathbf 1\{T_{a_0/2}<\textstyle\frac{t}{2}\}\right] +\mathbf P_{\tilde x}^M[T_{a_0/2}\geq \textstyle\frac{t}{2}]
\nonumber\\[1mm]
&=& \mathbf E_{\tilde x}^M\left\{\mathbf 1\{T_{a_0/2}<\textstyle\frac{t}{2}\}\mathbf E^M_{{X}^M_{T_{a_0/2}}}\left[{\rm e}^{-\int_0^{t-T_{a_0/2}} W(X_s^M,\tilde y)\,{\rm d}s}\right]\right\} +\mathbf P_{\tilde x}^M[T_{a_0/2}\geq \textstyle\frac{t}{2}]\nonumber\\[1mm]
&\leq & \mathbf E_{\tilde x}^M\left\{\mathbf 1\{T_{a_0/2}<\textstyle\frac{t}{2}\}\mathbf E_{{X}^M_{T_{a_0/2}}}^M\left[{\rm e}^{-\int_0^{{t}/{2}} W(X_s^M,\tilde y)\,{\rm d}s}\right]\right\} +\mathbf P_{\tilde x}^M[T_{a_0/2}\geq \textstyle\frac{t}{2}]\nonumber\\[1mm]
&\leq &
\mathbf E_{\tilde x}^M\left\{\mathbf 1\left\{T_{a_0/2}<\textstyle\frac{t}{2}\right\}\cdot
\sup_{\xi\in B(\tilde y, \frac{a_0}{2})}\mathbf E_\xi^M\left[{\rm e}^{-\int_0^{t/2}W(X_s^M,\tilde y)\,{\rm d}s}\right]\right\} +\mathbf P_{\tilde x}^M\left[T_{a_0/2}>\textstyle\frac{t}{2}\right].\nonumber\\
\label{eq:P3-1}
\end{eqnarray}

The expected value inside the supremum in the last line can be estimated as
\begin{eqnarray*}
\mathbf E_\xi^M\left[{\rm e}^{-\int_0^{t/2}W(X_s^M,\tilde y)\,{\rm d}s}\right] &=&
\mathbf E_\xi^M\left[{\rm e}^{-\int_0^{t/2}W(X_s^M,\tilde y)\,{\rm d}s} \mathbf 1\{\tau_{B(\tilde y, a_0)}\geq \textstyle\frac{t}{2}\}\right]\\
&&
+\mathbf E_\xi^M\left[{\rm e}^{-\int_0^{t/2}W(X_s^M,\tilde y)\,{\rm d}s} \mathbf 1\{\tau_{B(\tilde y, a_0)}<\textstyle\frac{t}{2}\}\right]\\
&\leq & {\rm e}^{-At/2}\mathbf P_\xi^M[\tau_{B(\tilde y, a_0)}\geq \textstyle\frac{t}{2}] + \mathbf P_\xi^M[\tau_{B(\tilde y, a_0)}<\textstyle\frac{t}{2}],
\end{eqnarray*}
where $\tau_{B(\tilde y, a_0)}$ denotes the first exit time of the process from the ball $B(\tilde y, a_0)$. This estimate is in the form $z\leq x(1-y)+y,$ with $x \in (0,1)$ and $y\in[0,1],$
so if we can prove that $y \leq p \in (0,1)$, then we will also have $z\leq x(1-p)+p.$

So far, all the estimates pertained to the projected stable process
$X_{(\gamma)}^M.$  Observe that for any open subset $F\subset \mathcal G_M,$  $x\in\mathcal G_M,$ and $s>0$  one has $\mathbf P_x^M[\tau_F<s] \leq \mathbf P_x[\tau_F<s],$ so that it is enough to
estimate $\mathbf P_\xi[\tau_{B(\tilde y, a_0)}<\frac{t}{2}]$ for the nonprojected process (i.e. the $\gamma-$stable process on $\mathcal G$).
As in present case  $\xi\in B(\tilde y, a_0/2),$ we have
\[\mathbf P_\xi[\tau_{B(\tilde y, a_0)}<\textstyle\frac{t}{2}]\leq
\mathbf P_\xi[\sup_{s\leq t/2} d(X_s,X_0)>\textstyle\frac{a_0}{2}],\]
which is not bigger than $\leq c {t}/{a_0^\gamma}$
in the jump stable ($\gamma<d_w$) case \cite[Lemma 4.3]{bib:BSS}, and not bigger that\linebreak $c{\rm e}^{-c_2(d_w) \left(\frac{a_0}{t^{1/d_w}}\right){d_w/(d_w-1)}}$ in the Brownian motion ($\gamma=d_w$) case \cite[Theorem 4.3]{bib:BP}, where $c=c(\gamma)$. In the sequel, we continue with the jump stable case only, the other case follows identically (easier in fact). The important feature of these estimates is that they do no longer depend on $M.$
Therefore, for any $0 < t \leq t_0'=\frac{1}{2}\,\frac{a_0^\gamma}{c},$
\[\sup_{\xi\in B(\tilde y, a_0/2)}\mathbf E_\xi^M\left[{\rm e}^{-\int_0^{t/2}W(X_s^M,\tilde y)\,{\rm d}s }\right] \leq
{\rm e}^{-At/2} \left(1-\frac{c t}{a_0^\gamma}\right)+ \frac{c t}{a_0^\gamma} \; \leq\; \frac{1}{2}\left({\rm e}^{-At/2}+1\right) .\]
Insert this estimate into (\ref{eq:P3-1}) and continue in the same vein. This time, use the observation that
\begin{eqnarray*}
\mathbf P^M_{\tilde x}[ T_{a_0/2}>\textstyle\frac{t}{2}] &\leq &
\mathbf P^M_{\tilde x}[X^M_{t/2}\notin B(\tilde y, a_0/2)] = 1 - \mathbf P^M_{\tilde x}[X^M_{t/2} \in B(\tilde y, a_0/2)] \leq 1- \mathbf P_{\tilde x}[X_{t/2} \in B(\tilde y, a_0/2)].
\end{eqnarray*}
Recall that the transition density $p$ of the (nonreflected) jump $\gamma-$stable process in $\cG$ enjoys the estimate (see \cite[Theorem 3.1]{bib:BSS} and \cite[Theorem 1.1]{bib:CKu})
\begin{equation}\label{eq:density}
\frac{1}{c^{(1)}}\, \min\left(\frac{t}{d(x,y)^{d+\gamma}},
  t^{-d/\gamma}\right) \leq  p(t,x,y)\leq c^{(1)}\,\min\left(\frac{t}{d(x,y)^{d+\gamma}},
  t^{-d/\gamma}\right),
  \end{equation}
  with certain positive constant $c^{(1)}=c^{(1)}(\gamma).$ Moreover, since $d(\tilde x, \tilde y)\leq b,$ for $z\in B(\tilde y, a_0/2)$ one has $d(\tilde x, z) \leq b+a_0/2.$ Consequently, from (\ref{eq:density}) we get, with $c^{(2)}=c^{(2)}(\gamma, d),$
  \[\mathbf P_{\tilde x}[X_{t/2}\in B(\tilde y, a_0/2)]\geq m(B(\tilde y, a_0/2))\cdot \inf _{z\in B(\tilde y, a_0/2)} p(t/2,\tilde x, z)\geq c^{(2)} a_0^d \min \left (\frac{t}{(b+\frac{a_0}{2})^{d+\gamma}}, t^{-d/\gamma}\right).\]
When $t \leq t_0''=(b+a_0/2)^\gamma,$ then this estimate is
  \[\mathbf P_{\tilde x}[X_{t/2}\in B(\tilde y, a_0/2)]\geq \frac{c^{(2)}ta_0^d }{(b+\frac{a_0}{2})^{d+\gamma}},\]
  being a number smaller than one if $t \leq t_0'''=\frac{(b+a_0/2)^{d+\gamma}}{c^{(2)} a_0^d}.$
  Collecting all of these, we obtain
  \begin{eqnarray*}
  \mathbf E_{\tilde x}^M\left[{\rm e}^{-\int_0^t W(X_s^M,\tilde y)\,{\rm d}s}\right] \leq \frac{1}{2}\, \frac{c^{(2)} ta_0^d}{(b+\frac{a_0}{2})^{d+\gamma}}\left({\rm e}^{-At/2}+1\right) +\left(1-\frac{c^{(2)}ta_0^d}{(b+\frac{a_0}{2})^{d+\gamma}}\right).
  \end{eqnarray*}
  For $t^*=\min(t_0',t_0'',t_0''')$
  we get an estimate in the form
\[\frac{p}{2} ({\rm e}^{-At^*/2} + 1) + (1-p)=:1-2c_{4.7},\]
with $p<1.$ To conclude, we choose $\tau_0=2t^*.$ The resulting constant $c_{4.7}$ is strictly positive and depends on $ A, a_0,b,\gamma$ only. The proof is complete. \hfill $\Box$

\subsection{The upper bound for general potentials}

We have the following  statement, matching  in general setting the lower bound
from Theorem \ref{th:lower}.

\begin{theo}\label{th:upper-gen}
Suppose $X$ is a subordinate Brownian motion in $\mathcal G$ via a complete subordinator $S$ with Laplace exponent $\phi$ of the form  (\ref{def:phi}), satisfying {\bf (U1)}, {\bf (U2)}, or {\bf (U3)}. Let the profile $W$ satisfy {\bf (W1)} -- {\bf (W4)}, and suppose that for certain $\theta>0$  there is a number $ K \in [0,\infty)$ such that
\[K= \liminf_{d(x,y)\to \infty} W(x,y) d(x,y)^{d+\theta}.\]
Let $\gamma=d_w$ (under {\bf (U1)}) or $\gamma=\alpha_1$ (under {\bf (U2)} and {\bf (U3)}).
Then there exist constants $E_1,E_1'>0$ such that:\\
\p (i) when $\gamma<\theta$ then
\[\limsup_{t\to\infty}\frac{\log L(t)}{t^{d/(d+\gamma)}}\leq -E_1 \nu^{\gamma/(d+\gamma)},\] \\
\p (ii) when $\gamma=\theta$ then
\[\limsup_{t\to\infty}\frac{\log L(t)}{t^{d/(d+\gamma)}}\leq -E_1 \nu^{\gamma/(d+\gamma)}-E_1'\nu,\]
\p (iii) when $\gamma >\theta$ then
\[\limsup_{t\to\infty}\frac{\log L(t)}{t^{d/(d+\theta)}}\leq -E_1' \nu.\]
\end{theo}

\noindent{\bf Proof.} The statements follows immediately from Proposition \ref{prop:long_range_2} and Theorem \ref{gorne-sbm}. Both statements are true in present setting, so we can use the arithmetic mean of both the bounds. The constants we
obtain are: $E_1=\frac{1}{2}D_1,$ $E_1'=\frac{K}{2c_{2.2} e^K}.$
\hfill$\Box$

\begin{rem} \label{re:func_upper}
\rm Identical statements hold true for the Feynman-Kac functionals. We skip the proof.
\end{rem}

\section{The Lifschitz tail for the integrated density of states}
In this section we transform the bounds from Theorems  \ref{th:lower} and  \ref{th:upper-gen} into bounds concerning the rate of decay of $l$ near zero. This is done by means of an exponential Tauberian theorem \cite[Th. 2.1]{bib:F}.
\begin{theo}\label{th:lif-upper}
Suppose that the assumptions of Theorem \ref{th:lower} are met. Then there exist constants $\widetilde C_1,...,\widetilde C_4>0$ such that:\\
\p (i) when $\beta<\theta$ then
\[\liminf_{x\to 0} x^{d/\beta}\log l([0,x]) \geq -\widetilde C_1 \nu,\]
\p (ii) when $\beta=\theta$ then
\[\liminf_{x\to 0} x^{d/\beta}\log l([0,x]) \geq -\widetilde C_2 \nu - \widetilde C_3 \nu^{1+d/\beta},\]
\p (iii) when $\beta>\theta$ then
\[\liminf_{x\to 0} x^{d/\theta}\log l([0,x]) \geq - \widetilde C_4 \nu^{1+d/\theta}.\]
\end{theo}

\begin{theo}\label{th:lif-lower}Suppose that the assumptions of Theorem \ref{th:upper-gen} are met. Then there exist constants $\widetilde D_1,...,\widetilde D_4>0$ such that:\\
\p (i) when $\gamma<\theta$ then
\[\limsup_{x\to 0} x^{d/\gamma}\log l([0,x])\leq -\widetilde D_1 \nu,\]
\p (ii) when $\gamma=\theta$ then
\[\limsup_{x\to 0} x^{d/\gamma}\log l([0,x]) \leq -\widetilde D_2 \nu - \widetilde D_3 \nu^{1+d/\gamma},\]
\p (iii) when $\gamma>\theta$ then
\[\limsup_{x\to 0} x^{d/\theta}\log l([0,x]) \leq - \widetilde D_4 \nu^{1+d/\theta}.\]
\end{theo}
The most interesting case is that of $\beta=\gamma,$ i.e. the case
when the lower and upper scaling exponents for $\phi$ coincide.
In this case, the rate of decay of $l([0,x])$ as $x \to 0^+$ is of order ${\rm e}^{-const\cdot x^{-d/\beta}}.$ Likewise, when $\theta <\beta$, i.e., when the behaviour of the potential at infinity
dominates the behaviour of the process, then the rate of decay of $l([0,x])$ is ${\rm e}^{-const\cdot x^{-d/\theta}}.$

\section{Examples}\label{sec:examples}

 At the very end, we give various examples of subordinators with Laplace exponents that are complete Bernstein functions satisfying the regularity assumptions needed for our work.

\begin{exam} {\rm We first discuss some examples of functions $\phi$ satisfying all of our assumptions.
\begin{itemize}
\item[(1)] \emph{Pure drift.} Let $\phi(\lambda)= b \lambda$, $b>0$. The corresponding subordinate process is just the Brownian motion with speed $b>0$. Clearly, the assumption \textbf{(U1)} is satisfied and \textbf{(L1)} holds with $\beta=d_w$.
\end{itemize}
\noindent
The next two examples are jump subordinators with drift.
\begin{itemize}
\item[(2)] \emph{Stable subordinator with drift.} Let $\phi(\lambda)=b\lambda + \lambda^{\gamma/d_w}$, $\gamma \in (0,d_w)$, $b>0$. Then the corresponding subordinator is a sum of a pure drift subordinator $b t$ and the pure jump $\gamma/d_w$-stable subordinator. In this case, \textbf{(L1)} and \textbf{(U2)} are satisfied with $\alpha_1=\alpha_2=\beta=\delta=\gamma$.
\item[(3)] Let $\phi(\lambda)=b\lambda + \lambda^{\gamma_1/d_w}[\log(1+\lambda)]^{\gamma_2/d_w}$, $\gamma_1 \in (0,d_w)$, $\gamma_2 \in (-\gamma_1,d_w-\gamma_1)$, $b>0$. In this case, we may take $\alpha_1=\beta=\gamma_1+\gamma_2$, $\alpha_2= \gamma_1$ and $\delta = (\gamma_1+d_w)/2$ in \textbf{(L1)} and \textbf{(U2)}.
\end{itemize}
\noindent
We now give examples of pure jump subordinators.
\begin{itemize}
\item[(4)] \emph{Mixture of purely jump stable subordinators.} Let $\phi(\lambda)=\sum_{i=1}^n \lambda^{\gamma_i/d_w}$, $\gamma_i \in (0,d_w)$, $n \in \N$. One can directly check that \textbf{(L1)} and \textbf{(U3)} hold with $\alpha_1 = \beta=\min_{i} \gamma_i$ and $\alpha_2=\delta=\max_{i} \gamma_i$.
\item[(5)] Let $\phi(\lambda)= (\lambda+\lambda^{\gamma_1/d_w})^{\gamma_2/d_w}$, $\gamma_1, \gamma_2 \in (0,d_w)$. The conditions \textbf{(L1)} and \textbf{(U3)} hold with $\alpha_1=\beta=(\gamma_1 \gamma_2)/d_w$ and $\alpha_2 = \delta=\gamma_2$.
\item[(6)] Let $\phi(\lambda)=\lambda^{\gamma_1/d_w}[\log(1+\lambda)]^{-\gamma_2/d_w}$, $\gamma_1 \in (0,d_w)$, $\gamma_2 \in (0,\gamma_1)$. One can check that both assumptions \textbf{(L1)} and \textbf{(U3)} are fulfilled for $\alpha_1=\alpha_2=\beta=\gamma_1-\gamma_2$  and $\delta=\gamma_1$.
\end{itemize}
\noindent
The last example satisfies our assumptions in part only. More precisely, it fulfils \textbf{(S1)}, \textbf{(S2)} and \textbf{(L1)}, but \textbf{(U1)}, \textbf{(U2)} nor \textbf{(U3)} do not hold.
\begin{itemize}
\item[(7)] \emph{Relativistic $\alpha/d_w$-stable subordinator.} Let $\phi(\lambda)=(\lambda+\vartheta^{d_w/\alpha})^{\alpha/d_w}-\vartheta$, $\alpha \in (0,d_w)$, $\vartheta>0$. The subordination via this subordinator leads to a very significant process called \emph{relativistic $\alpha$-stable}. Here \textbf{(L1)} holds with $\beta=d_w$, but neither of the conditions \textbf{(U1)}, \textbf{(U2)} nor \textbf{(U3)} is satisfied. Theorem \ref{th:lower} can be applied to this process with such $\beta$, but our Theorem \ref{gorne-sbm} does not cover this case. It can be conjectured that  appropriate upper bounds hold true with the same rate $\gamma=d_w$. However, proving this would require more specialized arguments customized to the specific properties of the relativistic stable process.
\end{itemize}
}
\end{exam}

We also provide examples of profile functions satisfying the assumptions of present paper.
\begin{exam}{\bf\cite[Example 4.1]{bib:KaPP}} {\rm
Fix $M_0 \in \mathbb Z$ and let the function $\psi:\cG_{M_0} \to [A,\infty)$ be such that $\psi \in L^1(\cG_{M_0},m),$ with $A>0.$  Define
$$
W(x,y):=
\left\{
\begin{array}{ll}
\psi(\pi_{M_0}(y)), & \mbox{when} \ x, y \in \Delta_{M_0}(z_0), \ \mbox{for some} \ z_0 \in \cG \backslash \mathcal V_{M_0}, \\
0,       & \mbox{otherwise} .
\end{array}
\right.
$$ It is established in \cite{bib:KaPP} that {\bf (W1)}--{\bf(W3)} hold true. Clearly, {\bf (W4)} holds as well.}
\end{exam}

\begin{exam}{\bf\cite[Example 4.2]{bib:KaPP}} {\rm
Let $\varphi:[0,\infty) \to [0,\infty)$ be a function satisfying the following conditions.
\begin{itemize}
\item[(1)] There exists $R>0$ such that $\varphi(x)=0$ for all $x \in (R,\infty)$.
\item[(2)] For every $y \in \cG$ one has $\varphi(d(\cdot,y))  \in \cK_{\loc}^X$.
\item[(3)] There exist numbers $a_0, A>0$ such that $\varphi(x)\geq A$ when $x<a_0.$
\end{itemize}
For such a function $\varphi$ we define
\begin{eqnarray}\label{ex:phi}
W(x,y):=\varphi(d(x,y)), \quad x, y \in \cG.
\end{eqnarray}
Again, {\bf (W1)}--{\bf(W3)} were verified in \cite{bib:KaPP},
and {\bf(W4)} is straightforward.
}\end{exam}

\medskip

 We also give an example of profile functions $W$ on $\cG \times \cG$ with unbounded support satisfying our assumptions. Such profiles can be realized as follows.

\begin{exam} {\rm
First we set  additional notation. Recall that for $M \in \mathbb Z_+$ and every $x \in (\cG \backslash \mathcal V_M) \cup \left\{0\right\}$ there is exactly one triangle (the so-called natural cell) of size $2^M$ in $\cG$, $\Delta_M(x)$, such that $x \in \Delta_M(x)$. If $x \in \mathcal V_M \backslash \left\{0\right\}$, then there are exactly two triangles $\Delta^{(1)}_M(x)$ and $\Delta^{(2)}_M(x)$ of size $2^M$ such that $\left\{x\right\} = \Delta^{(1)}_M(x) \cap \Delta^{(2)}_M(x)$. For every $x \in \cG$ we define
$$
r(x):= \sup_{\left\{p \in \mathcal V_0: \ d(x,y) \leq 1 \right\}} \sup \left\{M \in \mathbb Z_+: p \in \mathcal V_M\right\}.
$$
In particular, $r(0)=\infty$. One can check that for every $x \in \cG$ there is exactly one vertex $p_x \in \mathcal V_0$ such that $d(x,p_x) \leq 1$ and $r(x) = \sup\left\{M \in \mathbb Z_+: p_x \in \mathcal V_M\right\}$. For $x \in \cG$ and $M \in \mathbb Z_+$ we denote
$$
D_M(x)=\left\{
\begin{array}{ll}
\Delta^{(1)}_M(p_x) \cup \Delta^{(2)}_M(p_x) & \mbox{when} \quad M \leq r(x) < \infty, \\
\Delta_M(x) & \mbox{when} \quad M>r(x)  \ \mbox{or} \ p_x=0.
\end{array}
\right.
$$
Moreover, let $(a_n)_{n \in \mathbb Z_+}$ be a nonincreasing sequence of nonnegative numbers such that
\begin{equation} \label{eq:intexp}
\sum_{M=1}^{\infty} \sum_{n=[M/4]+1}^{\infty} 2^{nd} a_n < \infty.
\end{equation}
With the above notation we define
\begin{eqnarray*}
W(x,y) :=\left\{
\begin{array}{ll}
a_0 & \mbox{when} \quad x \in \cG, \ y \in D_0(x), \\
a_n & \mbox{when} \quad x \in \cG, \ y \in D_n(x) \backslash D_{n-1}(x),, \ n = 1, 2, 3, ... \ . \\
\end{array}
\right.
\end{eqnarray*}
Checking assumptions \textbf{(W1)} and \textbf{(W2)} for the profile $W$ is an easy exercise. The geometric condition \textbf{(W3)} can be established by  similar arguments as those in the justification in \cite[Example 4.3]{bib:KaPP}. To fulfil the condition \textbf{(W4)}, it is enough to assume that $a_0, a_1 >0$. The decay conditions for the profile $W$ as in (\ref{eq:cond-W}) and (\ref{eq:cond-W_up}) can be obtained by imposing some additional regularity on the elements of the sequence $(a_n)_{n \in \mathbb Z_+}$ for large $n$. For instance, by taking $a_{n} = 2^{- n (d+\theta)}$ with $0<\theta < \infty$, we immediately get
$$
2^{-d-\theta} =\liminf_{d(x,y) \to \infty} W(x,y) d(x,y)^{d+\theta} < \limsup_{d(x,y) \to \infty} W(x,y) d(x,y)^{d+\theta} = 1.
$$
By putting $a_n =0$ for $n \geq n_0 $, with some $2 \leq n_0 \in \mathbb Z_+$, we obtain a profile $W$ with bounded support.
}
\end{exam}

\section{Appendix: the enlargement of obstacles method for Markov processes with compact state-space}

The method of enlargement of obstacles was first introduced in
\cite{bib:Szn1} for diffusion processes on a compact state-space,
 evolving among killing obstacles. Its main
ingredient is an estimate comparing the principal eigenvalue of
the semigroup of such a process with the principal eigenvalue of
the process evolving in a modified environment -- with much bigger
obstacles. It has been proven that under appropriate conditions on
the process and on the configuration of the obstacle points, the
principal eigenvalue does not increase significantly after such a
modification, provided the principal eigenvalue of the initial
process was not too big. The method was generalized to some
non-diffusion Markov processes in \cite{bib:Kow}. We now need a version of these theorems for processes influenced by a killing potential with microscopic range, not microscopic killing
obstacles.

\subsection{The setting and the assumptions} \label{subsec:ass}
 Our initial setup consists of:
\begin{description}
\item[$*$] a compact linear metric space $({\cal T}, d)$ equipped with a doubling Radon measure $m,$ satisfying $m({\cal T})=1$.
    More precisely, we assume that there exist $r_0>0$ and $C_d\geq 1$ such that for any $x\in\mathcal T$ and $0<r<r_0$
    \begin{eqnarray}\label{doubling}  m(B(x,r))&\leq& C_d m(B(x, \frac{r}{3})),\end{eqnarray}
\item[$*$] a right-continuous, strong Markov process $X=\left(X_t, \pr_x\right)_{t \geq 0, \, x \in \cal T}$ on $\cal T$ with symmetric and strictly positive
transition density
 $p(t,x,y)$ with respect to $m$ such that $\forall\,\, t,$ $\int_{\cal T} p(t,x,x){\rm d}m(x)<\infty,$
\item[$*$] a potential profile $W:{\cal T}\times{\cal T}\to\mathbb{R}_+$ of finite range:
 a measurable function with support included in
$\{(x,y)\in\mathcal T\times\mathcal T: d(x,y)\leq a\epsilon\},$ where
$a>0, \epsilon>0$ are given, such that
\begin{equation} \label{eq:finiteW}
\mbox{for every} \quad t>0 \ \ \mbox{and} \ \  y \in {\cal T}, \quad \sup_{x \in {\cal T}} \ex_x \int_0^t W(X_s,y)\,{\rm  d}s < \infty.
\end{equation}

In applications, $a$ will be considered fixed and $\epsilon$ will tend to $0.$
\end{description}

Suppose $x_1,...,x_N\in \cal {T} $ are given points (`obstacles'),
then one defines the potential $V(x)$ as follows:
 \begin{equation}\label{vee} \mathcal{T}\ni x\mapsto
V(x)= \sum_{i=1}^N W(x,x_i). \end{equation} In applications, these
points will be random and  coming
 from a realisation of a Poisson point process $\cal N$ on $\mathcal{T}.$
Clearly, under the condition (\ref{eq:finiteW}) we have $\sup_{x \in {\cal T}} \ex_x \int_0^t V(X_s)\, {\rm d}s < \infty$, for every $t>0$.

Below we will study the process $X$ perturbed by the potential $V$. Formally, we consider the Feynman-Kac semigroup $(P^V_t)_{t \geq 0}$ on $L^2({\cal T}, m)$ consisting of symmetric operators
$$
P^V_t f(x) = \ex_x \left[e^{-\int_0^t V(X_s)\,{\rm d}s } f(X_t) \right], \quad f \in L^2({\cal T}, m), \ \ t > 0.
$$
Operators $P^V_t$ admit measurable, bounded and strictly positive kernels $p^V(t,x,y)$. Since also \linebreak $m({\cal T})=1<\infty$, all $P^V_t$ are of Hilbert-Schmidt type and have discrete spectra $\left\{e^{-\lambda_k(V) t}\right\}_{k=1}^{\infty}$, where $0 \leq \lambda_1(V) < \lambda_2(V) < ... \to \infty$ are eigenvalues of the generator of the semigroup $(P^V_t)_{t \geq 0}$. The corresponding eigenfunctions are denoted by $\varphi_k^V$. All $\lambda_k(V)$ have finite mutliplicity, the principal eigenvalue $\lambda_1(V)$ is simple, and the ground state eigenfunction $\varphi_1^V$  can be chosen to be strictly positive.

\smallskip

We intend to perform the following operation: for given  $b\gg a$ we would like to  replace the support
of the potential $V$ by a much larger set $\bigcup_{i=1}^N
\overline B(x_i,b\epsilon),$ and then to kill the initial process $X$ when it enters this bigger
set. Since $ \int_{\cal T} p(t,x,x)\,{\rm d}x<\infty,$ the semigroup of
this process again consists of symmetric Hilbert-Schmidt operators having discrete
spectrum. We are interested in comparing the smallest
eigenvalue of its generator with the principal eigenvalue $\lambda_1(V)$ of the process $X$ perturbed by the potential $V$. In general, we cannot enlarge every
obstacle  -- we need to restrict our attention
 to
those obstacles $x_i$ that are well-surrounded by other obstacles
(so-called {\em good obstacles}, see below).  Other obstacles will
be disregarded. Formally, we consider the sets
\begin{eqnarray}\label{thetabe}
\mathcal O_b =\bigcup_{x_i-{\rm good}} \overline{B}(x_i,b\epsilon), &&
\Theta_b=\mathcal T\setminus\mathcal{O}_b.
 \end{eqnarray}
 The process evolves now in the open set $\Theta_b$ and is killed when it enters $\mathcal O_b.$
Denote by $\lambda_1(b)$ the smallest eigenvalue of the generator
of this process.

The distinction between `good' and `bad' points is made as
follows.

\begin{defi}\label{defi:good-bad}
 Suppose $b,\delta,R >0$ are given, and let $x_1,...,x_N$ be given obstacle points. Then $x_{i_0}$ is called a good obstacle point if for all balls
$C=B(x_{i_0}, 10b\epsilon R^l)$ one has
\begin{equation}\label{good}
m\left(\bigcup_{i=1}^N \overline{B}(x_i,b\epsilon)\cap C\right)\geq
\frac{\delta}{C_d}\, m(C),
\end{equation}
($C_d$ is the constant from (\ref{doubling})) for all $l=0,1,2,...,$
as long as $10b\epsilon R^l<r_0.$ Otherwise, $x_{i_0}$ is called a bad obstacle point.
\end{defi}
Formally speaking, this notion depends on $b,\delta,R,$ but for
the time being we do not incorporate these parameters into the
notation.

 Balls with centers at bad obstacle points sum up to a
set with small volume.

\begin{lem}\emph{\cite[Lemma 1.3]{bib:Szn1}}\label{lembad}
\begin{equation}\label{badpoints}
m\big(\bigcup_{x_i-{\rm bad}}\overline{B}(x_i, b\epsilon)\big)\leq \delta.
\end{equation}
\end{lem}

\noindent We consider the following set of assumptions regarding
the process $X$ and the potential profile $W.$

\begin{description}
\item[(P1)]  There exists  $c_0>0 $ such that $ \sup_{x,y\in\mathcal T} p(1,x,y)\leq c_0.$
\end{description}

The remaining assumptions are concerned with recurrence properties
of the process. We require that for any fixed $a,b,$  $a\ll b$, $b\epsilon<r_0$
and  $\delta>0$ there exist constants $\tau_0, c_1, c_2,$ $
c_3,\alpha,\kappa>0,$ $R>3$ and a nonicreasing function
$\phi:(0,r_0)\to (0,1]$ such that:
\begin{description}
\item[(P2)] for $x,y\in\mathcal{T}$ with $d(x,y)\leq b\epsilon$ one
has
\[\mathbf{P}_x[\tau_{B(y,10 (R-2) b\epsilon)}<\frac{\tau_0\epsilon^\alpha}{2}]<c_1;\]
\item[(P3)]  when
$x,y\in\mathcal{T},$ and $d(x,y)\leq b\epsilon,$  then  $$
\mathbf{E}_x\big[{\rm e}^{-\int_0^{(\tau_0\epsilon^\alpha)/2} W( X_s,y)\,{\rm
d}s}\big]\leq 1-2c_1;$$
\item[(P4)] for $x,y\in\mathcal T$ satisfying $d(x,y)\leq  r\epsilon\leq
r_0$ one has
\[\mathbf{P}_x\big[T_{B(y, b\epsilon)}
\leq\frac{\tau_0\epsilon^\alpha}{2}]\geq \phi(r);\]
\item[(P5)] for
 $10 b\epsilon \leq \beta\leq\frac{r_0}{R},$
  any points $x,y\in\mathcal{T}$ with  $d(x,y)\leq\beta,$ and for
any compact subset $E\subset\mathcal{T}$ satisfying
$m(E\cap\overline{B}(y,\beta))\geq \delta/C_d \cdot
m(\overline{B}(y,\beta)) $ one has
\[\mathbf{P}_x[T_E<\tau_{B(y, R\beta)}]\geq c_2;\]
\item[(P6)] for $r<r_0/3,$ $A>3r$ and  $x,y\in\mathcal{T}$ satisfying
$d(x,y)\leq r$ one has
\[\mathbf{P}_x\left[X_{\tau_{B(y,r)}}\notin B(y, A)\right]\leq c_3\left(\frac{r}{A}\right)^{\kappa}.\]
\end{description}
  Assumption {\bf (P6)}  was first introduced in \cite{bib:Kow}.  It is typical for jump-type processes and was not needed
in the diffusion case.

\smallskip

As a preparatory step, we relate the expression involving the term
${\rm e}^{-\int _0^t V(X_s)\,{\rm d}s}$ to certain survival
probability of the process. As in \cite[page 171]{bib:Szn-env} (see also \cite{bib:kpp-point}), we
attach exponential clocks to each of the points $x_i,$ and then
kill the process $X_t$ once the quantity $A^i_t:=\int_0^t
W(X_s,x_i)\,{\rm d}s$ becomes bigger than this clock. More
precisely: given configuration of points $\left\{x_i\right\}_{i=1}^N$ we consider $N$ Poisson processes $(N^i(t))_{t \geq 0}$ with intensity $1$ on the probability spaces $(\Omega^i, \pr^i)$ and the product measure $\pr_z^N :=\pr_z \otimes (\otimes_{i=1}^N \pr^i)$ defined on the product of the canonical space for the process $X$ and spaces $\Omega^i$ of the processes $(N^i(t))_{t \geq 0}$, endowed with the product $\sigma-$algebra.
Product measures $\pr_z^N $  turn the canonical process $X$ and the canonical right
continuous counting processes $(N^i(t))_{t \geq 0}$ into independent processes,
distributed respectively as the initial process $X$ starting from $z$ and Poisson counting
processes with unit intensity, starting from $0$. They also satisfy the strong Markov property with respect to the appropriate filtration.

Let
\begin{eqnarray}\label{t-exponential}
T_i=\inf\{s\geq 0:\; N^i(A_s^i)\geq 1\}, &\mbox{and}&
T=\min_{i=1,...,N}T_i.
\end{eqnarray}
The following relation is central to our considerations:
\begin{eqnarray}\label{moments}
\pr^N_z[T>t\big|X] &= & \pr^N_z[\forall i=1,...,N\;\; T_i>t\big| X]
= \prod _i {\rm e}^{-\int_0^t W(X_s,x_i)\,{\rm d}s}=  {\rm
e}^{-\int_0^t V(X_s)\,{\rm d}s}.
\end{eqnarray}
Formally, in the next subsection we will be working with the measure $\pr^N_z$.
However, for simplicity, we do not indicate this in the notation, writing just $\pr_z$.

\subsection{The theorem comparing  the bottoms of the spectra}
The statement of the theorem together with its proof are very similar to that of \cite[Theorem 1.4]{bib:Szn1} and \cite[Theorem 1]{bib:Kow} -- the only difference is that we are concerned now with a killing
potential, slowing down the process, and not killing obstacles of small radius. To make the paper self-contained, we state the
theorem and briefly sketch the proof, indicating the changes that must be introduced and skipping the parts identical to those in
previous papers.

Before we state the  theorem, we make some
additional technical preparation, needed in the c\`adl\`ag case. For given  $K>\delta>0$  define
\begin{equation}\label{ce-em-delta}
C(K,\delta)= {\rm e}^K\left(1+c_0(1+\frac{K}{\delta})\right),
\end{equation}
where $c_0$ is the constant from the relation {\bf(P1)}. Suppose
that the number $R$ entering assumptions {\bf (P2), (P5)}
satisfies
\begin{equation}\label{erzero}
 \frac{c_3}{R^{\kappa}-1}\leq \frac{1}{8} \,C(K,\delta)^{-1}.
\end{equation} This can be done without loss of generality:
if {\bf (P2), (P5)} are satisfied with certain $R>0,$ then
they are satisfied for any $\tilde R>R.$

\begin{theo}\label{compare} Assume that the process $X_t$ is either:
\begin{enumerate}
\item[(i)] a diffusion satisfying {\bf (P1)} -- {\bf (P5)}
or
\item[(ii)] a discontinuous c\`adl\`ag process satisfying {\bf (P1)} -- {\bf (P6)}, with $R$ satisfying (\ref{erzero}).
    \end{enumerate}
    Let the numbers
$K>\delta>0,$ $ b\gg a$ be given. Then there exists $\epsilon_0=
\epsilon_0(a,b,\delta,K,c_0,c_1,c_2,c_3,\alpha,\kappa)$ \linebreak ($c_3$ and $\kappa$ not needed in the diffusion case) such that for any $\epsilon<\epsilon_0$ ($b \epsilon$ is
the  radius of obstacles in (\ref{thetabe})) one has
\begin{equation}\label{lambda}
\lambda_1(b)\wedge K\leq \lambda_1(V)\wedge
K+\delta.
\end{equation}
\end{theo}

\

\noindent{\bf Proof.} We prove (i) and (ii) simultaneously.

Let  $m_0$ be the smallest possible integer for which
\begin{eqnarray}\label{emzerozero}
(1-c_1c_2)^{m_0}\leq \frac{1}{8}\,C(K,\delta)^{-1} && \mbox{(diffusion case (i))}
\end{eqnarray}
or
\begin{eqnarray}\label{emzero}
(1-c_1c_2)^{\log_2m_0}\leq \frac{1}{8}\,C(K,\delta)^{-1} && \mbox{(non-diffusion case (ii))}.
\end{eqnarray}
Then we set \begin{equation}\label{dee} D:= 10b  R^{m_0}.
\end{equation}
 Denote $\lambda=\lambda_1(b)\wedge K-\delta.$ When $\lambda \leq 0,$ there is nothing to prove, so
we assume $\lambda>0.$
Our goal is to establish that
 \begin{equation}\label{glowne}\int_{\mathcal T} \mathbf E_x[{\rm e}^{\lambda T}] {\rm d} m(x)<\infty,
 \end{equation}
where $T$ is the stopping time introduced in (\ref{t-exponential}).

\smallskip

\noindent The proof of (\ref{glowne}) is divided in four steps.

\smallskip

 \noindent \textsc{Step 1.} For the stopping time
 \[ T_b=\inf\{t\geq 0: X_s\in \mathcal O_b\}\]
we just repeat  the argument from \cite[pages 231-232]{bib:Szn1}  to get
\begin{equation}\label{star}
\mathbf{E}_x[{\rm e}^{\lambda T_b}]\leq {\rm e}^{K}\left(1+ c_0
\left(1+\frac{K}{\delta}\right)\right)=C(K,\delta),
\end{equation}
for any $x\in\mathcal{T}.$
\medskip

\noindent \textsc{Step 2.}   Denote   $\sigma_D=\inf\{t\geq 0:
X_t\notin \bigcup_{i=1}^N B (x_i, D\epsilon)\}$ and $\tilde T=
T\wedge \sigma_D$.
In this step
 we
estimate $\mathbf{E}_x[{\rm e}^{\lambda \tilde T}]$, for $x \in \bigcup_{i=1}^N B (x_i, D\epsilon).$

Using  the Fubini theorem we write
\begin{equation}\label{pp}\mathbf{E}_x[{\rm e}^{\lambda \tilde T}] =
1+\int_0^\infty \lambda{\rm e}^{\lambda u}\,\mathbf{P}_x[\tilde
T>u]\;{\rm d} u, \quad x \in \bigcup_{i=1}^N B (x_i, D\epsilon).\end{equation} The event $\{\tilde T>u\}$ means
that up to time $u$ we have stayed inside $\bigcup_{i=1}^NB(x_i,D\epsilon)
$ and that $T$ did not happen up to this moment. If $T$ were to
occur before $u,$ we would have to enter the support of $V$ before
that time (an increase of some $A_t^i$ can happen only when the process
falls within the range of the potential $V$). Pick $x\in
\bigcup_{i=1}^N B(x_i,D\epsilon)$ and let $i_0$ be such that $x\in
B(x_{i_0}, D\epsilon).$ Then (assumption {\bf (P4)})
\[\mathbf{P}_x[T_{B(x_{i_0}, b\epsilon)}\leq \frac{\tau_0\epsilon^\alpha}{2}]
\geq \phi(D)\] and for $y\in B(x_{i_0},b\epsilon)$ (assumption {\bf (P3)})
\[\mathbf{P}_y [T\leq \frac{\tau_0\epsilon^\alpha}{2}] = \mathbf{E}_y\big[1-{\rm e}^{-\int_0^{\tau_0\epsilon^\alpha/{2}}\,V(X_s){\rm
d}s}\big]\geq 1- \mathbf{E}_y\big[{\rm
e}^{-\int_0^{\tau_0\epsilon^\alpha/2}W(X_s,x_{i_0})\,{\rm d}s }\big]\geq 2c_1.\]
Then from the strong Markov property it follows that, once $x \in \bigcup_{i=1}^N B (x_i, D\epsilon),$
\[\mathbf{P}_x[\tilde T\leq \tau_0\epsilon^\alpha]  \geq 2c_1\phi(D)\]
and from the ordinary Markov property we get
\[\mathbf{P}_x[\tilde T >u] \leq (1-2c_1\phi(D))^{[u/(\tau_0\epsilon^\alpha)]}\leq \frac{1}{1-2c_1\phi(D)}\,(1-2c_1\phi(D))^{u/(\tau_0\epsilon^\alpha)}, \quad x \in \bigcup_{i=1}^N B (x_i, D\epsilon)\]
(when $2c_1\phi(D)<1;$ when $2c_1\phi(D)=1$ then the quantity
estimated is equal to 0). We now insert this estimate into
(\ref{pp}) and proceed similarly as in \cite[page 231]{bib:Szn1}, obtaining that
for all $x \in \bigcup_{i=1}^N B (x_i, D\epsilon),$ as long as $K\tau_0\epsilon^\alpha < \log
(1-2c_1\phi(D))^{-1},$
\begin{eqnarray}\label{wazne1}
\mathbf{E}_x[{\rm e}^{\lambda \tilde T}]&\leq& 1+ \frac{K\tau_0\epsilon^\alpha}
{(1-2c_1\phi(D))[\log(1-2c_1\phi(D))^{-1}-K\tau_0\epsilon^\alpha]} \nonumber
\end{eqnarray}
Now we set
\[\epsilon_0=\frac{r_0}{4D}\wedge \inf\left\{\epsilon>0:\frac{K\tau_0\epsilon^\alpha}
{(1-2c_1\phi(D))[\log(1-2c_1\phi(D))^{-1}-K\tau_0\epsilon^\alpha]} > \frac{1}{8}\, C(K,\delta)^{-1}\right\},\]
so that when $\epsilon\leq \epsilon_0,$ we have
\begin{equation}\label{2102-1}
\mathbf{E}_x[{\rm e}^{\lambda \tilde T}]\leq 1+ \frac{1}{8}C(K,\delta)^{-1}\leq 2, \quad x \in \bigcup_{i=1}^N B (x_i, D\epsilon).
\end{equation}

\medskip

\noindent \textsc{Step 3.} In this step, we finally estimate
$\mathbf{E}_x[{\rm e}^{\lambda T}],$ using estimates (\ref{star})
and (\ref{wazne1}). Introduce the following sequence of stopping
times: $S_0=0$ and
 \begin{eqnarray*}
  S_1= T_b+  \tilde T\circ\theta_{T_b},
 && S_{n+1}= S_n+S_1\circ\theta_{S_n},\;\; n=1,2,... \, .
 \end{eqnarray*}
 Observe that each of the $S_n$'s (for $n=1,2,...$) is realized at the moment $t$ when either:\linebreak
 $X_{t}\in \mbox{supp}\,V$ (and the expression involving the potential gets bigger than the exponential clock)
  or $X_t\notin  \bigcup_{i=1}^N B(x_i, D\epsilon).$ Since $D>b>a,$
 these two possibilities are distinct. Therefore, it makes sense to define
 \[L=\inf\{n: X_{S_n}\in \mbox{supp}\,V\}\]
(with the convention $\inf\emptyset=\infty$). By modifying the argument leading to (\ref{koniec}) below combined with the Borel-Cantelli lemma, one can show that $L$ is finite $\pr_x$-a.s., for every $x \in \cal T$.
 We have $T\leq S_L,$ and therefore $\mathbf{E}_x[{\rm e}^{\lambda T}]\leq \mathbf{E}_x[{\rm e}^{\lambda S_L}].$ We can write:
 \begin{eqnarray*}
\mathbf{E}_x[{\rm e}^{\lambda S_L}] &\leq& 1+ \sum_{k=1}^\infty
\mathbf{E}_x[{\rm e}^{\lambda S_k}
 \mathbf{1}\{ L=k\} ]\\
 &=& 1+ \sum_{k=0}^\infty \mathbf{E}_x[{\rm e}^{\lambda S_{k+1}}
 \mathbf{1}\{L=k+1\} ] =(*).
\end{eqnarray*}
In view of the observation above, $L=l$ means that $X_{S_0},...,
X_{S_{l-1}}\in \left( \bigcup_{i\in I} B(x_i, D\epsilon)\right)^c $ and
$X_{S_l}\in \mbox{supp}\,V$. Therefore the estimate continues as
\begin{eqnarray*}
(*)&\leq & 1+\sum_{k=0}^\infty \mathbf{E}_x[\mathbf{E}_x[{\rm
e}^{\lambda
(S_k+S_1\circ\theta_{S_k})}\mathbf{1}\{X_{S_0},...,X_{S_k}\notin
\mbox{supp}\,V\}\big|\mathcal{F}_{S_k}]]\\
&=& 1+ \sum_{k=0} ^\infty \mathbf{E}_x\big[ {\rm e}^{\lambda
S_k}\mathbf{1}\{X_{S_0},...,X_{S_k}\notin \mbox{supp}\,V \}\cdot
\mathbf{E}_{X_{S_k}}[{\rm e}^{\lambda S_1}]\big].
 \end{eqnarray*}
For the last expectation, we first use the strong Markov property
and then inequalities  (\ref{star}), (\ref{wazne1}):
\begin{eqnarray*}
\mathbf{E}_{X_{S_k}}[{\rm e}^{\lambda
S_1}]&=&\mathbf{E}_{X_{S_k}}[{\rm e}^{\lambda(T_b+\tilde T\circ
\theta_{T_b})}] = \mathbf{E}_{X_{S_k}}[{\rm e}^{\lambda
T_b}\mathbf{E}_{X_{T_b}}[{\rm e}^{\lambda \tilde T}]]\\
&\leq & 2 C(K,\delta).
\end{eqnarray*} Consequently,
\[\mathbf{E}_x[{\rm e}^{\lambda T}]\leq  1+\sum_{k=0}^\infty \mathbf{E}_x\big[{\rm e}^{\lambda S_k}\mathbf{1}\{X_{S_0}\notin
\mbox{supp}\, V, ..., X_{S_k}\notin\mbox{supp}\,V\}\big]\cdot
2C(K,\delta).\]
 For $k\geq 0$ we denote
 \[a_k=\mathbf{E}_x\big[{\rm e}^{\lambda S_k}\mathbf{1}\{X_{S_0}\notin
\mbox{supp}\, V, ..., X_{S_k}\notin\mbox{supp}\,V\}]. \] Our goal
now is to show that for $k=1,2,...$ we have $a_k\leq \rho
a_{k-1},$ with certain constant  $\rho\in(0,1).$ This will do, as
then we will
  have \begin{equation}\label{koniec} \mathbf{E}_x[{\rm
e}^{\lambda T}] \leq 1+ \sum_{k=0}^\infty \rho^k \cdot 2
C(K,\delta)= 1+\frac{2C(K,\delta)}{1-\rho}<\infty.\end{equation}
 We can write, for $k\geq 1,$
\begin{eqnarray*}a_k&=&\mathbf{E}_x\big[{\rm e}^{\lambda
S_k}\mathbf{1}\{X_{S_0}\notin \mbox{supp}\, V, ...,
X_{S_k}\notin\mbox{supp}\,V\}\big]\phantom{blabla}\\&=&
\mathbf{E}_x \big[{\rm e}^{\lambda
S_{k-1}}\mathbf{1}\{X_{S_0}\notin \mbox{supp}\, V, ...,
X_{S_{k-1}}\notin\mbox{supp}\,V\}\mathbf{E}_{X_{S_{k-1}}}[{\rm e }
^{\lambda S_1}\mathbf{1}\{S_1\notin\mbox{supp}\,V\}
]\big].\end{eqnarray*} Therefore it suffices to find a universal
bound on $\mathbf{E}_z[{\rm e}^{\lambda
S_1}\mathbf{1}\{X_{S_1}\notin\mbox{supp}\,V\}],$ for $z\notin
\mbox{supp}\,V.$

Using again the strong Markov property,  the inequality $ab\leq
(a-1)+b$ (valid for $a\geq 1,$ $b\leq 1$), (\ref{star}), and
(\ref{wazne1}), we write, for $z\notin \mbox{supp}\,V:$
\begin{eqnarray}\label{kwadracik}
\mathbb{E}_{z}[{\rm e } ^{\lambda
S_1}\mathbf{1}\{S_1\notin\mbox{supp}\,V\}] &\leq&
\mathbf{E}_{z}[{\rm e}^{\lambda T_b }\mathbf{E}_{X_{T_b}}[{\rm
e}^{\lambda\tilde
T}\mathbf{1}\{X_{\tilde T}\notin\mbox{supp}\,V\}]]\nonumber\\
&\leq & \mathbf{E}_{z}\left[{\rm e}^{\lambda
T_b}\left(\mathbf{E}_{X_{T_b}}({\rm e}^{\lambda \tilde
T}-1)+\mathbf{P}_{X_{T_b}}(X_{\tilde T}\notin\mbox{supp}\,V)
\right)\right]\nonumber\\
 &\leq &
\frac{1}{8}+\mathbf{E}_{z}\left[
{\rm e}^{\lambda T_b}\mathbf{P}_{X_{T_b}}[X_{\tilde T}\notin\mbox{supp}\,V]\right]\nonumber\\
&\leq& \frac{1}{8}+ {C(K,\delta)} \sup_{x\in \mathcal {O}_b}
\mathbf{P}_x[X_{\tilde T}\notin \mbox{supp}\,V].
\end{eqnarray}
When $x\in\mathcal{O}_b,$ then $x$ lies at a distance at most $b\epsilon$
 from  a good point, say $x_{j_0}.$  We have:
 \begin{equation}\label{2201-1}
 \mathbf{P}_x[X_{\tilde T}\notin\mbox{supp}\,V]=\mathbf{P}_x[\tilde T=\sigma_D] =\mathbf{P}_x[T>\sigma_D]\leq \mathbf{P}_x[T>\tau_{B(x_{j_0},D\epsilon)}].
 \end{equation}
 Identically as in \cite[page 233]{bib:Szn1} and \cite[pages 742-743]{bib:Kow} we obtain that
 \begin{eqnarray}\label{2201-2}
 \forall \,k>0 &&\forall\, x\in\mathcal{T}\mbox{ such that } d(x,x_{j_0})\leq 10 b\epsilon R^l
 \nonumber\\
 & &\phantom{lll}\mbox{one has } \mathbf{P}_x[T<\tau_{B(x_{j_0}, 10 b\epsilon R^{l+1+k})}]\geq c_1c_2>0.
 \end{eqnarray}

\

In the  diffusive case, applying the strong Markov property at moments
\(\tau_{B(x_{j_0}. 10b \epsilon R^l)},\) \linebreak $l=1,2,...,m_0-1$ we get that
for $x\in B(x_{j_0},b\epsilon)$ one has
\[\mathbf P_x[T>\tau_{B(x_{j_0},D\epsilon)}]\leq (1-c_1c_2)^{m_0}\leq
\frac{1}{8}\, C(K,\delta)^{-1}.\]

In the non-diffusive case, we use estimates from \cite[pages 743-745]{bib:Kow} and get
\[\mathbf P_x[T>\tau_{B(x_{j_0},D\epsilon)}]\leq (1-c_1c_2)^{\log_2 m_0}+ \frac{c_3}{R^{\kappa}-1}\leq \frac{1}{4}\, C(K,\delta)^{-1}.\]

 In either case we have
\[\sup_{x\in\mathcal O_b}\mathbf{P}_x[X_{\tilde T}\notin \mbox{supp}\,V]\leq \frac{1}{4}\, C(K,\delta)^{-1},\]
which inserted into (\ref{kwadracik}) results in the estimate
\[\sup_{z\notin\mbox{\scriptsize supp}\,V} \mathbf{E}_{z}[{\rm e } ^{\lambda
S_1}\mathbf{1}\{S_1\notin\mbox{{supp}}\,V\}] \leq \frac{3}{8}\;(
=:\rho).\] Relation (\ref{koniec}) follows.

\

\noindent \textsc{Step 4. The conclusion.}  As the estimate obtained
in Step 3 is uniform in
$x\in\mathcal T,$ inequality (\ref{glowne}) follows as well. From the Fubini theorem
 we have
\[\int_{\cal T} \mathbf{E}_x[{\rm e}^{\lambda T}]\,{\rm d}m(x)= 1 +\int_0^\infty \lambda {\rm e}^{\lambda v}
\int_{\cal T}\mathbf{P}_x[T>v]\,{\rm d}m(x)\,{\rm d}v. \]
Hence,
\begin{eqnarray*}     \infty>
\int_0^\infty \lambda {\rm e}^{\lambda v} \int_{\cal
T}\mathbf{P}_x[T>v]\,{\rm d}m(x)\,{\rm d}v &\geq & \int_0^\infty
\lambda {\rm e}^{\lambda v} \sum_k \langle
\varphi^V_k,\mathbf{1}\rangle^2_{{L^2(\mathcal T,m)}} {\rm
e}^{-\lambda^V_k v}\,{\rm d}v\\
&\geq & \langle\varphi^V_1,\mathbf 1\rangle_{{L^2(\mathcal
T,m)}}^2\lambda \int_0^\infty {\rm e}^{(\lambda-\lambda^V_1)v}{\rm
d}v.
\end{eqnarray*}
Since $\langle\varphi^V_1,\mathbf 1\rangle_{{L^2(\mathcal
T,m)}}>0$, the last integral is finite, and so $\lambda<\lambda^V_1.$ The
proof is concluded.\hfill $\Box$

\bigskip

\noindent
\textbf{Acknowledgements.} The authors thank T. Grzywny for valuable discussions on subordinators.

\end{document}